\documentclass{article}

 \usepackage{natbib}
\usepackage{graphicx}
\usepackage[utf8]{inputenc}
\usepackage[a4paper,top=3cm,bottom=3cm,left=2cm,right=2cm,marginparwidth=1.75cm]{geometry}
\usepackage{hyperref}       
\usepackage{url}            
\usepackage{booktabs}       
\usepackage{amsfonts}       
\usepackage{amsmath}        
\usepackage{nicefrac}       
\usepackage{microtype}      
\usepackage{xcolor}         
\usepackage{verbatim}       
\usepackage{caption}
\usepackage{subcaption}
\usepackage[capitalise,noabbrev]{cleveref}

\newtheorem{thm}{Theorem}
\crefname{thm}{Theorem}{Theorems}
\crefname{lem}{Lemma}{Lemmas}
\crefname{cor}{Corollary}{Corollaries}
\crefname{prop}{Proposition}{Propositions}
\newtheorem{ass}{Assumption}
\crefname{ass}{Assumption}{Assumptions}
\crefname{ex}{Example}{Examples}
\crefname{rem}{Remark}{Remarks}
\crefname{equation}{}{}

\DeclareMathOperator*{\op}{op}

\DeclareMathOperator{\diag}{diag}

\DeclareMathOperator*{\as}{a.s.}

\newcommand{\E}{\mathbb{E}} 
\newcommand{\V}{\mathbb{V}} 
\newcommand{\R}{\mathbb{R}} 
\newcommand{\N}{\mathbb{N}} 

\begin{document}


\title{On Adaptive Stochastic Optimization for Streaming Data: \\ A Newton's Method with $\mathcal{O}(dN)$ Operations}

\author{  Antoine Godichon-Baggioni$^{1}$ and Nicklas Werge$^{2}$}

\date{}
 
\maketitle

\noindent $^{1}$ Laboratoire de Probabilit\'es, Statistique et Mod\'elisation,
       Sorbonne Universit\'e, Paris, France,\\ antoine.godichon\_baggioni@sorbonne-universite.fr \\

\noindent $^{2}$  Department of Mathematics and Computer Science,
       University of Southern Denmark,
       Odense, Denmark, werge@sdu.dk

\begin{abstract}
Stochastic optimization methods encounter new challenges in the realm of streaming, characterized by a continuous flow of large, high-dimensional data.
While first-order methods, like stochastic gradient descent, are the natural choice, they often struggle with ill-conditioned problems. In contrast, second-order methods, such as Newton's methods, offer a potential solution, but their computational demands render them impractical.
This paper introduces adaptive stochastic optimization methods that bridge the gap between addressing ill-conditioned problems while functioning in a streaming context.
Notably, we present an adaptive inversion-free Newton's method with a computational complexity matching that of first-order methods, $\mathcal{O}(dN)$, where $d$ represents the number of dimensions/features, and $N$ the number of data.
Theoretical analysis confirms their asymptotic efficiency, and empirical evidence demonstrates their effectiveness, especially in scenarios involving complex covariance structures and challenging initializations.
In particular, our adaptive Newton's methods outperform existing methods, while maintaining favorable computational efficiency.
\end{abstract}

\noindent\textbf{Keywords: } stochastic optimization, adaptive methods, Newton's method, online learning, large-scale, streaming

\section{Introduction} \label{sec::intro}

The focus of this paper is on the stochastic optimization problem, where the objective is to minimize a convex function $F:\R^{d}\rightarrow\R$ with $d\in\N$. The problem is formulated as follows:
\begin{equation}
\label{eq::so_problem}
\min_{\theta\in\R^{d}}\{F(\theta):=\E_{\xi\sim\Xi}[f(\theta;\xi)]\},
\end{equation}
where $f$ is a loss function, $\xi$ is a random variable following an unknown distribution $\Xi$, and $\theta$ is the parameter of interest.
The challenge in \cref{eq::so_problem} is widespread in machine learning applications \citep{kushner2003stochastic,shapiro2021lectures,bottou2018optimization,sutton2018reinforcement}. For instance, in the context of an input-output pair $\xi=(x,y)$, the function $f$ typically takes the form $f(\theta;\xi)=f(\theta;x,y)=l(h_{\theta}(x);y)$, where $l$ is a loss function onto $\R$ and $h_{\theta}$ is a prediction model parameterized by $\theta$.

\medskip
We address the stochastic optimization problem \cref{eq::so_problem} within a streaming context, where data are both large in size and dimensionality. Similar to prior work by \citet{godichon2023non,godichon2023learning}, streaming data continuously arrives in blocks, resembling time-varying mini-batches, as independent and identically distributed (i.i.d.) samples of the random variable $\xi$. More formally, we consider an endless sequence of i.i.d. copies: $\{\xi_{1,1},\ldots,\xi_{1,n_{1}}\}, \ldots,\{\xi_{t,1},\ldots,\xi_{t,n_{t}}\},\ldots$, where $\{\xi_{t,1},\ldots,\xi_{t,n_{t}}\}$ represents a block of $n_{t}$ data points arriving at time $t$. This setup mirrors the incremental and block-based nature of real-world streaming data.

\medskip

Our adaptive stochastic optimization methods go beyond the conventional stochastic gradient-based methods by incorporating a Hessian matrix approximation $A_{t}$ at each step $t$ to refine the descent direction. In a general form, these adaptive methods can be expressed recursively as:
\begin{equation} \label{eq::asom}
\theta_{t+1} = \theta_{t} - \gamma_{t+1} A_{t} \nabla_{\theta}f(\theta_{t};\xi_{t+1}), \quad \theta_{0}\in\R^{d},
\end{equation}
where $\nabla_{\theta}f(\theta_{t};\xi_{t+1})=n_{t+1}^{-1}\sum_{i=1}^{n_{t+1}} \nabla_{\theta}f(\theta_{t};\xi_{t+1,i})$, $(\nabla_{\theta}f(\theta_{t};\xi_{t+1,i}))$ is unbiased gradients in $\R^{d}$, $(A_{t})$ is a sequence of random matrices in $\R^{d\times d}$, and $(\gamma_{t})$ is the learning rate.

\medskip

Specially, if we set $A_{t}=\mathbb{I}_{d}$ and $n_{t}=1$, the update in \cref{eq::asom} reduces to the classical Robbins-Monro method \citep{robbins1951}, commonly known as Stochastic Gradient Descent (SGD). When $n_{t}\in\N$ (with $A_{t}=\mathbb{I}_{d}$), we obtain a streaming version of SGD, akin to time-varying mini-batch SGD, as considered in \citet{godichon2023non,godichon2023learning}. For Adagrad \citep{duchi2011adaptive}, $A_{t}$ serves as an estimate of the inverse square root of the diagonal of the variance of the gradients. Furthermore, the update in \cref{eq::asom} transforms into Newton’s method, when $A_{t}$ serves as an approximation of the inverse Hessian matrix $\nabla_{\theta}^{2}F(\theta_{t})$.

\medskip

The central question in this paper is twofold: \emph{Can we construct a sequence of Hessian approximations $(A_{t})$ in a manner that is both computationally efficient and ensures the robustness of our adaptive methods to ill-conditioned problems?}

\paragraph{Contributions.} In this work, we present adaptive stochastic optimization methods capable of robustly handling ill-conditioned problems while ensuring computational efficiency in streaming contexts. These adaptive methods dynamically adjust learning per-dimension, leveraging historical gradient and Hessian information. Additionally, we propose iterative weighted average versions of our adaptive methods. These acceleration techniques both provide variance-reduction during learning and accelerated convergence. Theoretical analysis establishes their asymptotic efficiencies, encompassing strong consistency, rate of convergence, and asymptotic normality. Empirical evidence further validates their effectiveness, particularly in scenarios with complex covariance structures and challenging initializations.

\medskip

A noteworthy contribution of our work is the introduction of inversion-free adaptive Newton's methods, designed to match the computational complexity of first-order methods---$\mathcal{O}(dN_{t})$, where $N_{t}=\sum_{i=1}^{t}n_{i}$ is the total quantity of data up to time $t$. These adaptive Newton's methods not only achieve the computational efficiency of first-order methods but also incorporate acceleration techniques for enhanced convergence, while harnessing the power of second-order information.

\paragraph{Related work.} Stochastic optimization and adaptive methods have been extensively researched, as evident in works such as \citet{bottou2018optimization,chau2022inversion}. Theoretical investigations into SGD span topics from in-depth non-asymptotic analysis to its asymptotic efficiency \citep{moulines2011non,kushner2003stochastic,toulis2017asymptotic,pelletier1998almost,fabian1968asymptotic,pelletier2000asymptotic,gadat2022asymptotic,nemirovski2009robust,lacoste2012simpler}. A noteworthy extension of SGD is the concept of averaging, known for its role in accelerating convergence. This averaging scheme, referred to as Polyak-Ruppert averaging or averaged SGD (ASGD), was introduced by \citet{ruppert1988efficient,PolyakJud92}. They demonstrated that using a learning rate with slower decays, combined with uniform averaging, robustly leads to information-theoretically optimal asymptotic variance. While these estimates are known to be asymptotically efficient \citep{pelletier2000asymptotic}, their non-asymptotic properties have been thoroughly investigated \citep{moulines2011non,needell2014stochastic,gadat2023optimal}. However, it's important to note that this averaging concept can be sensitive to ill-conditioned problems among others, leading to sub-optimal performance in practice \citep{leluc2023asymptotic,boyer2023asymptotic}.

\medskip

To address this practical challenge, recent strategies have emerged to enhance the performance of stochastic optimization methods, focusing on adaptive approaches. These methods involve tuning the learning rate, also known as the step-size sequence, through strategies that adapt to the gradient. One of the most well-known adaptive techniques is Adagrad \citep{duchi2011adaptive}, which incorporates an estimation of the square root of the inverse of the gradient's covariance into the step-size. Subsequently, this method has undergone various modifications and improvements. Notable among these adaptations are RMSProp \citep{tieleman2012rmsprop}, ADAM \citep{kingma2015adam}, AdaDelta \citep{zeiler2012adadelta}, NADAM \citep{dozat2016incorporating}, and AMSGrad \citep{reddi2018on}. Nevertheless, these adaptive methods do not fully tackle the challenge of poor conditioning. Another limitation of these methods is their reliance on information solely from the diagonal of the gradient covariance estimator. Consequently, in scenarios with strong correlations, this restricted information may result in sub-optimal practical outcomes.

\medskip

To address this, an alternative approach involves considering inversion-free stochastic Newton's methods \citep{bercu2020efficient,boyer2023asymptotic,leluc2023asymptotic}, where an estimate of the inverse of the Hessian is integrated into the step-size. Alternatively, stochastic Gauss-Newton methods \citep{cenac2020efficient,bercu2021stochastic} can be employed. These stochastic Newton's methods, relying on the Sherman-Morrison formula \citep{sherman1950adjustment},\footnote{Sherman-Morrison's formula is also known as Riccati’s equation for matrix inversion \citep{duflo2013random}.} require a specific form of the Hessian. Nevertheless, they find applications in various scenarios, including linear, logistic, softmax, and ridge regressions \citep{bercu2020efficient,boyer2023asymptotic,godichon2024recursive}, as well as tasks such as estimation of the geometric median \citep{GW2023}, non-linear regression \citep{cenac2020efficient}, and optimal transport \citep{bercu2021stochastic}.

\medskip

Our adaptive stochastic optimization methods aim to integrate the strengths of acceleration techniques (weighted Polyak-Ruppert averaging), adaptive methods, and stochastic Newton's methods. We believe that this integration provides an effective solution to solving the challenges posed by ill-conditioned problems in a streaming context.

\paragraph{Organization.} This paper is organized as follows: \Cref{sec::framework} presents the underlying theoretical framework within we analyse our adaptive stochastic optimization methods. In \cref{sec::asom}, we analyse our adaptive stochastic optimization methods and its weighted averaged version in \cref{sec::asom::wa}. In
\cref{sec::app}, we apply our adaptive methods to Adagrad and Newton's method. In particular, \cref{sec::app} details the construction of our adaptive Newton's methods with $\mathcal{O}(dN_{t})$ operations. In \cref{sec::exp}, we present our experimental results, demonstrating the efficiency of our proposed methods. 

\paragraph{Notations.} We represent the Euclidean norm as $\lVert\cdot\rVert$ and the operator norm as $\lVert\cdot\rVert_{\op}$. The notation $M\succ0$ indicates that $M$ is positive definite, while $M\succeq0$ indicates that it is positive semi-definite. The minimum and maximum eigenvalues of matrix $M$ are denoted by $\lambda_{\min}(M)$ and $\lambda_{\max}(M)$, respectively.

\section{Underlying Theoretical Framework} \label{sec::framework}
In this section, we provide the theoretical framework that underpins our analysis. Our objective is to solve the stochastic optimization problem in \cref{eq::so_problem}, while operating within a streaming contexts. As a reminder, we consider stochastic optimization problems of the form:
\begin{equation*}
\min_{\theta\in\R^{d}}\{F(\theta):=\E_{\xi\sim\Xi}[f(\theta;\xi)]\},
\end{equation*}
where $\xi$ is a random variable sampled following an unknown distribution $\Xi$.

\medskip

To lay the foundation for our analyses, we introduce three key assumptions. These assumptions, contingent upon the differentiability of the function $F$, serve as the bedrock for our theoretical framework. These assumptions are standard in the realms of stochastic optimization, stochastic approximation, and adaptive methods \citep{bottou2018optimization,leluc2023asymptotic,boyer2023asymptotic,kushner2003stochastic,godichon2019online,godichon2019lp,benveniste-book90,duflo2013random,GBT2023}.
\begin{ass}\label{ass::1}
For almost any $\xi$, the function $f(\cdot;\xi)$ is differentiable and there exists non-negative constants $C$ and $C'$ for all $\theta\in\R^{d}$, such that
\begin{equation} \label{ass::1::eq}
\E[\lVert\nabla_{\theta}f (\theta;\xi)\rVert^{2}]\leq C+C'( F(\theta)-F(\theta^{*})).
\end{equation}
In addition, there exists $\theta^{*}\in\R^{d}$ such that $\nabla_{\theta} F(\theta^{*})=0$, and the functional $\Sigma : \theta \rightarrow \E[\nabla_{\theta}f(\theta;\xi)\nabla_{\theta}f (\theta;\xi)^{\top}]$ is continuous at $\theta^{*}$. 
\end{ass}

\medskip

In \cref{ass::1}, we do not confine $\E[\lVert\nabla_{\theta}f(\theta;\xi)\rVert^{2}]$ by a constant or the squared errors $\lVert\theta-\theta^{}\rVert^{2}$. Instead, we utilize the functional error $F(\theta)-F(\theta^{*})$,  a condition known as expected smoothness \citep{gower2019sgd,gazagnadou2019optimal,gower2021stochastic}. Moreover, when $C=0$, \cref{ass::1::eq} is known as the weak growth condition \citep{vaswani2019fast,nguyen2018sgd}. Notably, it is worth mentioning that, in the context of $\mu$-strong convexity of the functional $F$, the squared errors condition implies the functional error condition, as $\lVert\theta-\theta^{*}\rVert^{2}\leq\nicefrac{2}{\mu}(F(\theta)-F(\theta^{*}))$ for any $\theta\in\R^{d}$.

\medskip

In order to ensure the strong consistency of our method's estimates, we invoke a second key assumption. This assumption allows the use of a second-order Taylor expansion of the functional $G$ and is based on the following hypothesis:
\begin{ass}\label{ass::2}
The functional $F$ is twice-continuously differentiable with uniformly bounded Hessian, i.e., there exists $L_{\nabla F}$ such that $\lVert\nabla^{2}_{\theta}F(\theta)\rVert_{\op} \leq L_{\nabla F}$.
\end{ass}
Note that this implies, among other things, that the gradient of $F$ is $L_{\nabla F}$-Lipschitz.
The third assumption pertains to the uniqueness of the minimizer $\theta^{*}$ of the functional $F$.
\begin{ass}\label{ass::3}
The functional $F$ is locally strongly convex; $\lambda_{\min} := \lambda_{\min}(\nabla_{\theta}^{2}F(\theta^{*})) > 0$.
\end{ass}

\section{Adaptive Stochastic Optimization Methods} \label{sec::asom}
For clarity, our main discussion revolves around constant mini-batches of size $n$ (instead of time-varying mini-batches $n_{t}$). This approach enables us to intricately address the nuances of the streaming data setting while upholding the conceptual clarity of our core findings. However, it's crucial to emphasize that we offer translations and adaptations of these discussions for the scenario of time-varying mini-batches in \cref{sec::appendix}. The motivation for considering time-varying mini-batches stems from recent work by \citet{godichon2023learning,godichon2023non}, which demonstrated that increasing mini-batches can accelerate convergence and break long- and short-term dependence structures.

\medskip

Throughout the paper, we consider constant mini-batches of size $n$, i.e., at each time $t$, $n$ i.i.d copies of $\xi$ denoted by $\xi_{t}=\{\xi_{t,1},\ldots,\xi_{t,n}\}$ arrives. Our adaptive stochastic optimization methods, as defined in \cref{eq::asom}, can recursively be written as
\begin{equation*} 
\theta_{t+1} = \theta_{t} - \gamma_{t+1} A_{t}\nabla_{\theta}f(\theta_{t};\xi_{t+1}), \quad \theta_{0}\in\R^{d},
\end{equation*}
where $\nabla_{\theta}f(\theta_{t};\xi_{t+1})=n^{-1}\sum_{i=1}^{n} \nabla_{\theta}f(\theta_{t};\xi_{t+1,i})$. We assume the construction of a filtration $(\mathcal{F}_{t})$ such that $\theta_{t}$ and $A_{t}$ are $\mathcal{F}_{t}$-measurable, and $\xi_{t+1}=\{\xi_{t+1,1},\ldots,\xi_{t+1,n}\}$ are independent from $\mathcal{F}_{t}$. Let $N_{t}$ denote the total number of data processed up to time $t$, i.e., $N_{t}=nt$.

\medskip

Our goal is to recursively update $\theta_{t}$ at each time step $t$ to integrate the most recent information. For the subsequent discussion, we assume that the learning rate $(\gamma_{t})$ and the sequence of random matrices $(A_{t})$ satisfy the following conditions:
\begin{align}\label{cond::step}
\sum_{t \geq 1}\gamma_{t}\lambda_{\min}(A_{t-1})=+\infty \; \as, \quad \text{and} \quad \sum_{t \geq 1}\gamma_{t}^{2}\lambda_{\max}(A_{t-1})^{2}<+\infty \; \as
\end{align}
In \cref{sec::app}, we will delve into the modifications required in the methods to fulfill these conditions. Additional insights can be found in works such as \cite{boyer2023asymptotic,GBT2023}. In all the sequel, we take $\gamma_{t}=C_{\gamma}t^{-\gamma}$ with $C_{\gamma}>0$ and $\gamma\in(1/2,1)$ for the sake of simplicity. However, one can also take $\gamma_{t} = C_{\gamma}(t+t_{0})^{-\gamma}$ with $t_{0}\in\N$, and all the theoretical results remain true.

\medskip

The following theorem establishes the strong consistency of our adaptive stochastic gradient estimates $(\theta_{t})$ derived from \cref{eq::asom}.
\begin{thm} \label{theo::ps}
Suppose \cref{ass::1,ass::2,ass::3} hold, along with the conditions in \cref{cond::step}. Then, $\theta_{t}$ converges almost surely to $\theta^{*}$.
\end{thm}

To ascertain the rate of convergence of our adaptive stochastic gradient estimates $(\theta_{t})$, we assume that the sequence of random matrices $A_{t}$ converges to $A\succ0$.
\begin{ass}\label{ass::4}
The random matrix $A_{t}$ converges almost surely to a positive definite matrix $A$.
\end{ass}

For instance, in Newton's methods, the matrix $A$ represents the inverse Hessian, and in the case of Adagrad, it corresponds to the inverse of the square root of the diagonal of the gradient's variance. Note that once \cref{theo::ps} is fulfilled, the strong consistency of $\theta_{t}$ often implies the consistency of $A_{t}$. 

\begin{thm} \label{theo::rate}
Suppose \cref{ass::1,ass::2,ass::3,ass::4} hold, along with the conditions in \cref{cond::step}. In addition, assume there exist positive constants $C_{\eta}$ and $\eta>\frac{1}{\gamma}-1$ such that for all $\theta\in\R^{d}$,
\begin{equation}\label{eq::momenteta}
\mathbb{E}\left[\lVert \nabla_{\theta}f(\theta;\xi) \rVert^{2+2\eta} \right] \leq C_{\eta} \left( 1+ F(\theta) - F(\theta^{*}) \right)^{1+\eta}.
\end{equation}
Then, 
\begin{equation*}
\lVert \theta_{t} - \theta^{*} \rVert^{2} = \mathcal{O} \left( \frac{\ln(N_{t})}{N_{t}^{\gamma}} \right) \; \as
\end{equation*}
\end{thm}

\section{The Weighted Averaged Version} \label{sec::asom::wa}
The weighted averaged version of our adaptive stochastic optimization methods in \cref{eq::asom} is defined for $w\geq0$ as follows:
\begin{equation} \label{eq::asom-wa}
\theta_{t,w}= \frac{1}{\sum_{i=0}^{t-1}\ln(i+1)^{w}} \sum_{i=0}^{t-1}\ln(i+1)^{w}\theta_{i},
\end{equation}
which can be written recursively as
\begin{equation*}
\theta_{t+1,w} = \theta_{t,w} + \frac{\ln(t+1)^{w}}{\sum_{i=0}^{t}\ln(i+1)^{w}}(\theta_{t}-\theta_{t,w}).
\end{equation*}
This weighted averaging in \cref{eq::asom-wa} enhances the optimization by adaptively assigning more weight to the latest estimates of $(\theta_{t})$. The logarithmic weighting emphasizes recent estimates, presumed to be more accurate, while providing robustness against sub-optimal initializations \citep{mokkadem2011generalization,boyer2023asymptotic}. Observe that taking $w=0$ leads to the usual Polyak-Ruppert averaging scheme \citep{ruppert1988efficient,PolyakJud92,godichon2023non}.

\medskip

To establish the convergence rate of the weighted averaged version of our adaptive stochastic optimization methods in \cref{eq::asom}, we begin by introducing a new assumption.
\begin{ass}\label{ass::5}
There exists positive constants $L_{\eta}$ and $\eta$ such that for all $\theta\in\mathcal{B}(\theta^{*},\eta)$,
\begin{equation*}
\lVert \nabla_{\theta}F(\theta) - \nabla_{\theta}^{2}F(\theta^{*})(\theta - \theta^{*}) \rVert \leq  L_{\eta} \lVert \theta - \theta^{*} \rVert^{2}.
\end{equation*}
\end{ass}
\Cref{ass::5} is satisfied as soon as the Hessian of $F$ is locally Lipschitz on a neighborhood around $\theta^{*}$. Coupled with \cref{ass::2}, this imply there exists a positive constant $L_{\delta}$ such that for all $\theta\in\R^{d}$,
\begin{equation*}
\lVert \nabla_{\theta}F(\theta) - \nabla_{\theta}^{2}F(\theta^{*})(\theta-\theta^{*}) \rVert \leq  L_{\delta} \lVert \theta - \theta^{*} \rVert^{2}.
\end{equation*}

The following result establish the rate of convergence and the optimal asymptotic normality of the weighted averaged estimates $(\theta_{t,w})$.
\begin{thm} \label{theotlc}
Suppose \cref{ass::1,ass::2,ass::3,ass::4,ass::5} hold, along with inequality \cref{eq::momenteta}. In addition, assume there exists a positive constant $\nu$ such that
\begin{equation} \label{eq::A_rate_convergence}
\lVert A_{t} - A \rVert_{\op} = \mathcal{O} \left( \frac{1}{t^{\nu}} \right) \; \as
\end{equation}
Then,
\begin{equation*}
 \lVert \theta_{t,w} - \theta^{*} \rVert^{2} = \left\lbrace \begin{array}{lr} \mathcal{O} \left( \frac{\ln(N_{t})}{N_{t}^{\gamma+2\nu}} \right) \; \as &\text{if} \; 2\nu+\gamma \leq 1, \\
\mathcal{O} \left( \frac{\ln(N_{t})}{N_{t}} \right) \;  \as &\text{if} \; 2\nu +  \gamma > 1.
\end{array} \right. 
\end{equation*}
Moreover, if $2\nu+\gamma>1$, then
\begin{equation*}
\sqrt{N_{t}}(\theta_{t,w} - \theta^{*}) \xrightarrow[t\to+\infty]{\mathcal{L}} \mathcal{N}(0,\nabla_{\theta}^{2}F(\theta^{*})^{-1}\Sigma\nabla_{\theta}^{2}F(\theta^{*})^{-1}).
\end{equation*}
\end{thm}
To establish strong results, such as the asymptotic efficiency of the weighted average estimates $(\theta_{t,w})$, the sequence of random matrices $A_{t}$ should exhibit a (weak) rate of convergence, as outlined in \cref{eq::A_rate_convergence}. In simpler terms, achieving a satisfactory rate of convergence of $A_{t}$ leads to the asymptotic efficiency of the weighted average estimates $(\theta_{t,w})$.

\medskip

However, to establish asymptotic efficiency without relying on a (weak) rate of convergence of $A_{t}$, one can also consider the following theorem:
\begin{thm}\label{theowithoutrate}
Suppose \cref{ass::1,ass::2,ass::3,ass::4,ass::5} hold, along with inequality \eqref{eq::momenteta}. In addition, assume there exists a positive constant $v'>1/2$ such that
\begin{equation}\label{equalitystrange}
\frac{1}{\sum_{i=0}^{t-1}\ln(i+1)^{w}} \sum_{i=0}^{t-1} \ln(i+1)^{w+1/2+\delta} \lVert A_{i+1}^{-1}-A_{i}^{-1} \rVert_{\op}  (i+1)^{\frac{\gamma}{2}} = \mathcal{O} \left(\frac{1}{t^{v'}} \right) \; \as,
\end{equation}
for some $\delta>0$. Then
\begin{equation*}
\lVert \theta_{t,w} - \theta^{*} \rVert^{2} = \mathcal{O} \left( \frac{\ln(N_{t})}{N_{t}}\right) \; \as \quad \text{and} \quad \sqrt{N_{t}}(\theta_{t,w}-\theta^{*}) \xrightarrow[t\to+\infty]{\mathcal{L}} \mathcal{N}(0,\nabla_{\theta}^{2}F(\theta^{*})^{-1}\Sigma \nabla_{\theta}^{2}F(\theta^{*})^{-1}).
\end{equation*}
\end{thm}

Note that, while condition \cref{equalitystrange} may appear unusual, it is straightforward to verify in practice. The proofs of \cref{theo::adagrad,theo::newton::nd} provide insights into practical methods for checking this condition.

\section{Applications to  Newton's Method} \label{sec::app}
In this section, we apply our adaptive stochastic optimization methodology, as detailed in \cref{sec::asom,sec::asom::wa}, to (stochastic) Newton's methods. Here, we present inversion-free adaptive Newton's methods explicitly designed to align with the computational complexity of first-order optimization methods. Specifically, we present a weighted average inversion-free adaptive Newton method that seamlessly integrates the strengths of both approaches. Additionally, it's worth noting that we propose a novel streaming variant of Adagrad, along with its weighted average counterpart, in \cref{sec::app::adagrad}.

\medskip

To overcome the computational challenges linked to Hessian inversion, we propose a variant of the stochastic Newton's method that entirely avoids Hessian inversion. In \cref{sec::app::newton:direct}, our adaptive stochastic optimization methodology is applied to the stochastic Newton's method, resulting in a direct streaming stochastic Newton's method requiring $\mathcal{O}(d^{2}N_{t})$ operations. Next, in \cref{sec::app::newton:nd}, we introduce a weighted Hessian estimate demanding only $\mathcal{O}(dN_{t})$ operations. Finally, in \cref{sec::app::newton:nd:wa}, we accelerate this stochastic Newton's method through weighted Polyak-Ruppert averaging.

\subsection{Direct Streaming Stochastic Newton's Method} \label{sec::app::newton:direct}
In the special case of stochastic Newton's methods, one can obtain the asymptotic efficiency without averaging by taking a step sequence of the form $\gamma_{t}=\nicefrac{1}{t}$.\footnote{Observe, in the increasing batch-size case in \cref{sec::appendix}, one should take $\gamma_{t} = \nicefrac{n_{t}}{N_{t}}$.} The streaming stochastic Newton algorithm is defined by the update:
\begin{equation} \label{eq::newton}
\theta_{t+1} = \theta_{t} - \frac{1}{t+1} \bar{H}_{t}^{-1} \nabla_{\theta}f(\theta_{t};\xi_{t+1}), \quad \theta_{0}\in\R^{d},
\end{equation}
where $\nabla_{\theta}f(\theta_{t};\xi_{t+1})=n^{-1}\sum_{i=1}^{n} \nabla_{\theta}f(\theta_{t};\xi_{t+1,i})$ and $\bar{H}_{t}$ is an approximation of the Hessian of $F$. Specifically, we consider Hessian estimates $\bar{H}_{t}=N_{t}^{-1}H_{t}$ of the form
\begin{equation*}
H_{t} = H_{0} + \sum_{i=1}^{t}\sum_{j=1}^{n} \alpha_{i,j} \Phi_{i,j}\Phi_{i,j}^{\top},
\end{equation*}
with $H_{0}$ symmetric and positive, $\alpha_{i,j}= \alpha(\theta_{i-1};\xi_{i,j})$, and $\Phi_{i,j}=\Phi(\theta_{i-1};\xi_{i,j})$. A computationally-efficient estimate of $H_{t}$ inverse can be derived using the Riccati/Sherman-Morrisson’s formula \citep{duflo2013random,sherman1950adjustment} used $n$ times, i.e for all $j = 1 , \ldots , n$,
\begin{equation*}
H_{t-1+ \frac{j}{n}}^{-1} = H_{t-1 +  \frac{j-1}{n}}^{-1} -   \alpha_{t,j} \left( 1 + \alpha_{t,j} \Phi_{t,j}^{\top} H_{t-1 +  \frac{j-1}{n}}^{-1}  \Phi_{t,j}\right)^{-1} H_{t-1 +  \frac{j-1}{n}}^{-1}  \Phi_{t,j} \Phi_{t,j}^{\top}H_{t-1 +  \frac{j-1}{n}}^{-1} .
\end{equation*}
with the convention $H_{t -1-  \frac{ 1}{n}}^{-1} = H_{t-1}$.
The explicit construction of the recursive estimates of the inverse Hessian is detailed in various applications, including linear, logistic, softmax, and ridge regressions \citep{bercu2020efficient,boyer2023asymptotic,godichon2024recursive}. Additionally, these methods are applied to tasks such as the estimation of the geometric median \citep{GW2023}, non-linear regression \citep{cenac2020efficient}, and optimal transport \citep{bercu2021stochastic}.

\medskip

The asymptotic efficiency of the streaming version of the stochastic Newton's method can now be articulated as follows:
\begin{thm} \label{theo::newton::direct}
Suppose \cref{ass::1,ass::2,,ass::3,ass::5} hold, along with inequality \cref{eq::momenteta}. Then, $\theta_{t}$ converges almost surely to $\theta^{*}$.
In addition, assume $\bar{H}_{t}^{-1}$ converges almost surely to $\nabla_{\theta}^{2}F(\theta^{*})^{-1}$. Then,
\begin{equation*}
\lVert \theta_{t} - \theta^{*} \rVert^{2} = \mathcal{O} \left( \frac{\ln(N_{t})}{N_{t}} \right) \; \as
\end{equation*}
Moreover, assume there exists a positive constant $\nu$ such that $\lVert \bar{H}_{t}^{-1} - \nabla_{\theta}^{2}F(\theta^{*})^{-1} \rVert_{\op} = \mathcal{O}(\nicefrac{1}{t^{\nu}})\; \as$. Then
\begin{equation*}
\sqrt{N_{t}}(\theta_{t}-\theta^{*}) \xrightarrow[t \to + \infty]{\mathcal{L}} \mathcal{N} (0,\nabla_{\theta}^{2}F(\theta^{*})^{-1}\Sigma\nabla_{\theta}^{2}F(\theta^{*})^{-1}).
\end{equation*}
\end{thm}

\subsection{Streaming Stochastic Newton's methods with possibly $\mathcal{O}(dN_{t})$ operations} \label{sec::app::newton:nd}
The direct stochastic Newton's method presented in \cref{sec::app::newton:direct} is associated with computational costs of $\mathcal{O}(d^{2}N_{t})$, which can be computationally expensive, especially in high-dimensional streaming settings. To address this challenge, we introduce the streaming stochastic Newton's method using weighted Hessian estimates:
\begin{equation} \label{eq::newton::nd}
\theta_{t+1} = \theta_{t} - \frac{1}{t+1} \bar{H}_{t,w'}^{-1} \nabla_{\theta}f(\theta_{t};\xi_{t+1}), \quad \theta_{0}\in\R^{d},
\end{equation}
where $\nabla_{\theta}f(\theta_{t};\xi_{t+1})=n^{-1}\sum_{i=1}^{n} \nabla_{\theta}f(\theta_{t};\xi_{t+1,i})$ and $\bar{H}_{t,w'}=N_{t,Z}^{-1}H_{t,w'}$ with
\begin{equation} \label{eq::ssn:nd:hessian}
H_{t,w'} = H_{0,w'} + \sum_{i=1}^{t}\ln (i+1)^{w'}\sum_{j=1}^{n} Z_{i,j} \left(  \iota_{i,j} \tilde{e}_{i,j}\tilde{e}_{i,j}^{\top}  +   \alpha_{i,j} \Phi_{i,j}\Phi_{i,j}^{\top} \right),
\end{equation}
with $N_{t,Z} = 1+\sum_{i=1}^{t}\ln(i+1)^{w'}\sum_{j=1}^{n} Z_{i,j}$,  $H_{0}$ symmetric and positive, $w' \geq 0$, and $Z_{i,j}$ are i.i.d with $Z_{i,j} \sim \mathcal{B}(p)$ for some $p\in (0,1]$.
In addition, let $N_{t,k,Z} = (1+ \sum_{i=1}^{t-1}\sum_{j=1}^{n}Z_{i,j} + \sum_{j=1}^{k} Z_{t,j})$, $\iota_{i,j} = c_{\iota}  N_{i,j,Z}^{-\iota}$ for $\iota \in (0, 1/2  )$, and $e_{i,j} $ be the ($N_{i,j,Z}$ modulo $d+1$)-th component of the canonical basis. Observe that the term $\iota_{t}$ enables to control the smallest eigenvalue of $\bar{H}_{t,w'}$ while $\ln (t+1)^{w'}$ enables us to give more weights to the lastest updates $\alpha_{t,j}\Phi_{t,j}\Phi_{t,j}^{\top}$, which are supposed to be better since $(\theta_{t})$ converges to $\theta^{*}$. The add of the random variables $Z_{i,j}$ enables us to play with the computational cost. As explain later, taking $p=1$ leads to a natural recursive estimate of the Hessian in \cref{sec::app::newton:direct}.

\medskip

We now discuss about $Z_{i,j}$. Let us recall that with the help of Riccati's formula \citep{duflo2013random}, one can update the inverse of $H_{t+1,w'}$ as follows:
\begin{equation*}
\left\lbrace \begin{array}{lll}
 {H}_{t+\frac{1}{n},w'}^{-1}  &=   {H}_{t,w'}^{-1} - \frac{ Z_{t+1,1} \iota_{t+1,1} }{ 1+ \iota_{t+1,1}e_{t+1,1}\tilde{H}_{t,w'}^{-1} e_{t+1,1}^{T} }  {H}_{t,w'}^{-1}e_{t+1,1}e_{t+1,1}^{T} {H}_{t,w'}^{-1}     \\
  &\vdots \\
 {H}_{t + 1,w'}^{-1}  &=   {H}_{t+ \frac{n-1}{n},w'}^{-1} - \frac{ Z_{t+1,n} \iota_{t+1,n} }{ 1+ \iota_{t+1,n}e_{t+1,n} {H}_{t,w'}^{-1} e_{t+1,n}^{T} }  {H}_{t,w'}^{-1}e_{t+1,n}e_{t+1,n}^{T} {H}_{t,w'}^{-1}     \\
 {H}_{t+\frac{1}{n},w'}^{-1} &= {H}_{t+1,w'}^{-1} - \frac{Z_{t+1,1} \ln (t+1)^{w'}\alpha_{t+1,1}}{ 1+ \ln (t+1)^{w'}\alpha_{t+1,1} \Phi_{t+1,1}^{T} {H}_{t+1,w'}^{-1}\Phi_{t+1,1}}  {H}_{t+1,w'}^{-1} \Phi_{t+1,1}\Phi_{t+1,1}^{T} {H}_{t+1,w'}^{-1} \\
 {H}_{t+\frac{2}{n},w'}^{-1} &= {H}_{t+\frac{1}{n},w'}^{-1} - \frac{Z_{t+1,2}\ln (t+1)^{w'}\alpha_{t+1,2}}{ 1+ \ln (t+1)^{w'}\alpha_{t+1,2} \Phi_{t+1,2}^{T} {H}_{t+\frac{1}{n},w'}^{-1}\Phi_{t+1,2} }  {H}_{t+\frac{1}{n},w'}^{-1} \Phi_{t+1,2}\Phi_{t+1,2}^{T} {H}_{t+\frac{1}{n},w'}^{-1} \\
 &\vdots   \\
 {H}_{t+1,w'}^{-1} &=  {H}_{t+\frac{n-1}{n},w'}^{-1} - \frac{Z_{t+1,n}\ln (t+1)^{w'}\alpha_{t+1,n}}{ 1+ \ln (t+1)^{w'}\alpha_{t+1,n} \Phi_{t+1,n}^{T} {H}_{t+\frac{n-1}{n},w'}^{-1}\Phi_{t+1,n} }  {H}_{t+\frac{n-1}{n},w'}^{-1} \Phi_{t+1,n}\Phi_{t+1,n}^{T} {H}_{t+\frac{n-1}{n},w'}^{-1}.
\end{array} \right.
\end{equation*}
Then, the update of $H_{t+1}^{-1}$ only costs, in average, $\mathcal{O}(pd^{2}n)$ operations leading to a total number of operations of order (in average);
\begin{equation*}
  \underbrace{pd^{2}N_{t}}_{\text{estimating the inverse Hessian}} + \underbrace{dN_{t}}_{\text{estimating the gradient}} + \underbrace{\frac{d^{2}N_{t}}{n} }_{\text{multiplication of Hessian and gradient estimates}}.
\end{equation*}
Thus, one can play with the value of $p$ to reduce the cost of the update of the inverse of the Hessian. Indeed, one can obtain an averaged computational cost at time $t$   of order $\mathcal{O}( dN_{t} )$ operations  taking $p=d^{-1}$ and $n=d$. In other words, it is possible to obtain a stochastic Newton's method with only $\mathcal{O}(dN_{t})$ operations, which still is asymptotically efficient. {In all the sequel, for the sake of simplicity, we suppose that for any $\theta \in \mathbb{R}^{d}$,
\begin{equation}\label{hessegale}
\nabla^{2}_{\theta}F(\theta) = \mathbb{E}\left[ \alpha (\theta ; \xi ) \Phi(\theta ; \xi)\Phi(\theta ; \xi)^{T} \right] .
\end{equation}
}

\begin{thm}\label{theo::newton::nd}
Suppose \cref{ass::1,ass::2,ass::3,ass::5} hold, along with inequalities \cref{eq::momenteta,hessegale}.  {In addition, assume that for almost any $\xi$, there exists positive constants $C_{\eta'}$ and $\eta'>1$ such that for all $\theta\in\R^{d}$,
\begin{equation*}
\E[\lVert\alpha(\theta;\xi)\Phi(\theta;\xi)\Phi(\theta;\xi)^{\top}\rVert^{\eta'}] \leq C_{\eta'}^{\eta'}.
\end{equation*}}
Then, 
\begin{equation*}
\lVert \theta_{t} - \theta^{*} \rVert^{2} = \mathcal{O} \left( \frac{\ln(N_{t})}{N_{t}} \right) \; \as 
\end{equation*}
Moreover, suppose that the Hessian of $F$ is locally $L_{\nabla^{2}F}$-Lipschitz on a neighborhood around $\theta^{*}$ and that $\eta ' \geq 2$. Then
\begin{equation*}
\sqrt{N_{t}}(\theta_{t}-\theta^{*}) \xrightarrow[t \to + \infty]{\mathcal{L}} \mathcal{N}(0,\nabla_{\theta}^{2}F(\theta^{*})^{-1}\Sigma\nabla_{\theta}^{2}F(\theta^{*})^{-1}).
\end{equation*}
\end{thm}
Observe that contrary to Theorem \ref{theotlc}, no restriction on $\nu$ is necessary.
\subsection{Weighted Averaged Version of Streaming Stochastic Newton's methods with possibly $\mathcal{O}(dN_{t})$ operations} \label{sec::app::newton:nd:wa}
Although the streaming Newton's methods is very performant, it can be quite sensitive to bad initialization since the learning rate can be too small \citep{boyer2023asymptotic}. In order to overcome this, one can consider the weighted averaged version, given by
\begin{align}
\theta_{t+1} =& \theta_{t} -\gamma_{t+1} \bar{S}_{t,w'}^{-1} \nabla_{\theta}f(\theta_{t};\xi_{t+1}), \; \theta_{0}\in\R^{d}, \label{def::SN} \\
\theta_{t+1,w} =& \theta_{t,w} + \frac{\ln(t+1)^{w}}{\sum_{i=0}^{t}\ln(i+1)^{w}} (\theta_{t}-\theta_{t,w}), \label{def::WASN}
\end{align}
where $\nabla_{\theta}f(\theta_{t};\xi_{t+1})=n^{-1}\sum_{i=1}^{n} \nabla_{\theta}f(\theta_{t};\xi_{t+1,i})$, $\gamma_{t} = C_{\gamma}t^{-\gamma}$, and $\bar{S}_{t,w'}=N_{t,Z}^{-1}S_{t,w'}$ with
\begin{equation*}
S_{t,w'} =  S_{0,w'} + \sum_{i=1}^{t}\ln(i+1)^{w'} \sum_{j=1}^{n} Z_{i,j}\left( \iota_{i}  e_{i,j}e_{i,j}^{\top} +  \alpha_{i,j} \Phi_{i,j}\Phi_{i,j}^{\top} \right),
\end{equation*}
with $S_{0w'}$ symmetric and positive, $N_{t,Z}=1+\sum_{i=1}^{t}\ln(i+1)^{w'}\sum_{j=1}^{n} Z_{i,j}$, $\iota_{i} = c_{\iota}  N_{i,j,Z}^{-\iota}$, and $\iota \in (0,\gamma-1/2)$. One can follow  the same recursive scheme as for $H_{t,w'}^{-1}$ to update the inverse of $S_{t,w'}^{-1}$. Indeed, the only difference between $S_{t,w'}^{-1}$ and $H_{t,w'}^{-1}$ is the choice of the estimate of $\theta^{*}$.

\begin{thm} \label{theo::newton::wasn}
Suppose \cref{ass::1,ass::2,ass::3,ass::5} hold, along with inequalites \cref{eq::momenteta,hessegale}.  {In addition, assume that for almost any $\xi$, there exists positive constants $C_{\eta'}$ and $\eta'>1$ such that for all $\theta\in\R^{d}$,
\begin{equation*}
\E[\lVert\alpha(\theta;\xi)\Phi(\theta;\xi)\Phi(\theta;\xi)^{\top}\rVert^{\eta'}] \leq C_{\eta'}^{\eta'}.
\end{equation*}}
Then, 
\begin{equation*}
\lVert \theta_{t} - \theta^{*} \rVert^{2} = \mathcal{O} \left( \frac{\ln(N_{t})}{N_{t}^{\gamma }} \right) \; \as, \quad \quad \lVert \theta_{t,w} - \theta^{*} \rVert^{2} = \mathcal{O} \left( \frac{\ln(N_{t})}{N_{t}}\right) \; \as,
\end{equation*}
and
\begin{equation*}
\sqrt{N_{t}}(\theta_{t,w}-\theta^{*}) \xrightarrow[t \to + \infty]{\mathcal{L}}\mathcal{N}(0,\nabla_{\theta}^{2}F(\theta^{*})^{-1}\Sigma\nabla_{\theta}^{2}F(\theta^{*})^{-1}).
\end{equation*}
\end{thm}

\section{Experiments} \label{sec::exp}

In this section, we empirically evaluate our adaptive stochastic optimization methods, focusing on two fundamental problems in statistical optimization: least-squares regression and logistic regression.
For least-squares regression, our data points are represented as $\xi=(x,y)\in \mathbb{R}^{d}\times\mathbb{R}$, and we employ the functional $F(\theta) = \frac{1}{2}\mathbb{E}[(y-x^{\top}\theta)^{2}]$. In the case of logistic regression, our data points are $\xi=(x,y)\in\mathbb{R}^{d}\times\{-1,1\}$, and the corresponding functional is $F(\theta) = \mathbb{E}[\ln(1+\exp(x^{\top}\theta))-yx^{\top}\theta]$.

\medskip

To introduce complex covariance structures into our datasets, we adopt the experimental framework detailed in \citet{boyer2023asymptotic}. This involves modeling the covariance of our feature vector $x$ as follows:
\begin{equation*}
x \sim \mathcal{N}\left(0, M\diag\left(\frac{i^{2}}{d^{2}}\right)_{i=1,\dots,d}M^{\top}\right).
\end{equation*}
Here, $M$ represents a random orthogonal matrix. This choice of covariate distribution, influenced by the action of $M$, allows us to introduce strong correlations between the coordinates of $x$. This variation in data structure enables us to assess the adaptability of our method under diverse conditions.

\medskip

For our experiments, we deliberately set $d=100$ to emphasize the challenges posed by high dimensionality. This choice serves to highlight the scalability and robustness of our proposed methods when handling large-dimensional datasets. Within this setting, the Hessian associated with the model exhibits a wide range of eigenvalues, with the largest eigenvalue being a thousand times larger than the smallest one.

\medskip

The weight parameters for both the estimates and Hessian approximations are set to $2$, i.e., $w=w'=2$. While higher values of $w$ and $w'$ would enhance adaptability, the chosen setting serves as a proof-of-concept.

\medskip

We set the mini-batch size $n$ equal to the dimension $d$, which has the implication that our adaptive Newton's method, incorporating second-order information, updates the Hessian less frequently than the first-order SGD and Adagrad algorithms.\footnote{The streaming variant of Adagrad, along with its weighted average counterpart, is detailed in \cref{sec::app::adagrad}.} Surprisingly, our results clearly demonstrate that despite fewer Hessian updates, the Newton's method perform exceptionally well. Notably, when dealing with highly correlated data, the Adagrad algorithm's adaptive step-size becomes less effective, whereas the Newton's methods excel. Particularly, in scenarios involving less-than-ideal initializations (as depicted on the right side of the figures), both Newton's methods demonstrate outstanding performance.

\medskip

Leveraging this setup, we demonstrate the adaptability of our methods when dealing with high-dimensional datasets featuring complex covariance structures. Our experiments underscore the efficiency of our adaptive Newton's method in comparison to first-order gradient methods and highlight our methods' state-of-the-art performance in terms of both convergence speed and accuracy when compared to existing methods.

\medskip

In our applications, we investigate various optimization methods, encompassing SGD, Adagrad, our streaming Adagrad detailed in 
\cref{sec::app::adagrad}, along with their weighted Polyak-Ruppert averages. Additionally, we explore our streaming stochastic Newton (SSN) from \cref{sec::app::newton:nd} and our weighted Polyak-Ruppert averaged streaming stochastic Newton (WASSN) from \cref{sec::app::newton:nd:wa}; for SNN and WASSN, we set $p=1$ (i.e., $\mathcal{O}(d^{2}N_{t})$ computations) and $p=1/d$ (i.e., $\mathcal{O}(dN_{t})$ computations) to explore the loss of not updating the hole Hessian at each step.

\subsection{Least-Squares Regression} \label{sec::exp::lin_reg}

In the context of least-squares regression, we aim to evaluate the performance of our adaptive stochastic optimization methods for fitting linear models to the data. This entails modeling a linear relationship where the dependent variable $y$ is expressed as a linear combination of the feature vector $x$ and a parameter vector $\theta^{*}$. We adopt $\theta^{*}$ as our target parameter vector, specifically defined as $\theta^{*} = (-d/2,\dots, d/2)^{\top}$ \citep{boyer2023asymptotic}.

\medskip

\begin{figure}[!t]
    \centering
    \includegraphics[width=0.49\linewidth]{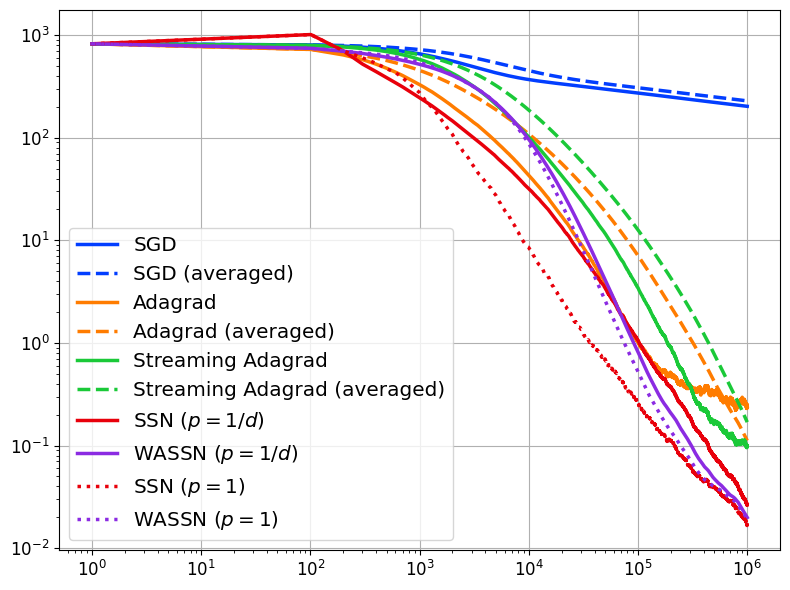}
    \includegraphics[width=0.49\linewidth]{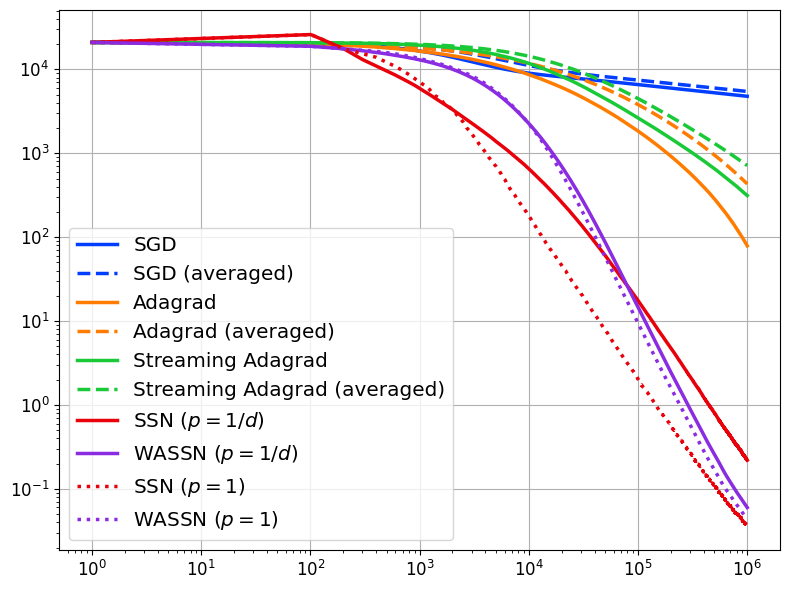}
     \caption{Least-Squares Regression: Mean-squared error of the distance to the optimum $\theta^{*}$, plotted against the sample size of $1.000.000$, for various initializations. The initial points $\theta_{0}$ are generated as $\theta_{0}=\theta^{*}(1+rU)$, where $U$ is a uniform random variable on the unit sphere of $\R^{d}$, and $r$ takes values of $1$ (left) or $5$ (right). Each curve reports $\lVert\theta_{t}-\theta^{*}\rVert$ averaged over $50$ different epochs, with a different initial point drawn for each sample.}
    \label{fig:lr}
\end{figure}

In \cref{fig:lr}, we present the quadratic mean error of different estimates, considering two types of initializations. Notably, we observe that both the Adagrad and Newton's methods exhibit faster convergence rates when compared to the standard Stochastic Gradient Descent (SGD). This enhanced performance is attributed to their innate capability to manage the diagonal structure of the Hessian matrix, which comprises eigenvalues at different scales. It's important to note that while the Adagrad algorithm adapts its step size, it may be less effective when confronted with highly correlated data. Intriguingly, for scenarios involving less-than-ideal initializations (as depicted on the right side of the figure), both Newton's methods demonstrate outstanding performance.

\subsection{Logistic Regression} \label{sec::exp::log_reg}

In logistic regression, our focus shifts to evaluating the performance of our adaptive stochastic optimization methods within the realm of binary classification. Logistic regression models the probability of a data point belonging to one of two classes based on predictor variables. We utilize a sigmoid function to transform a linear combination of the feature vector $x$ and the parameter vector $\theta^{*}$ into class probabilities. Consistent with \citet{boyer2023asymptotic}, we choose $\theta^{*}\in\R^{d}$ with all components equal to one. Unlike the linear regression setting, logistic regression exhibits intrinsic non-linearity, which makes the impact of the covariance structures less clear.

\medskip

\begin{figure*}[!t]
    \centering
    \includegraphics[width=0.49\linewidth]{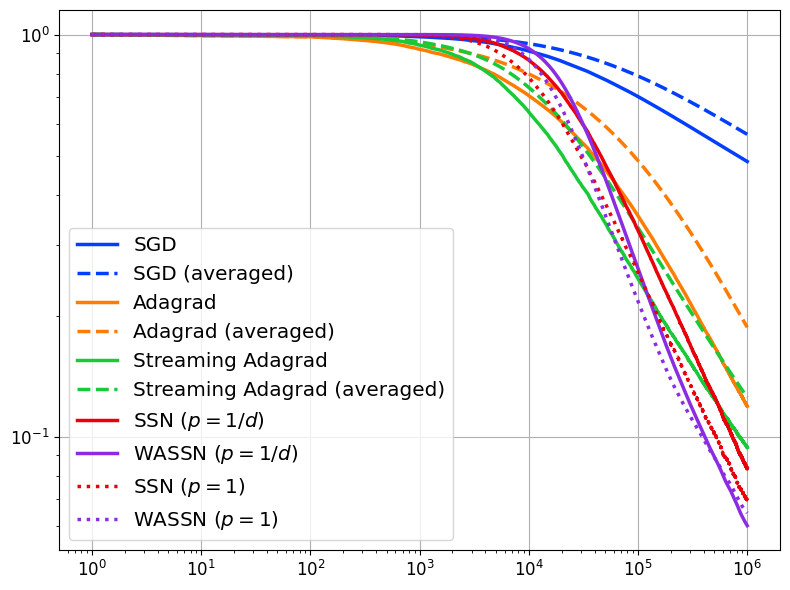}
    \includegraphics[width=0.49\linewidth]{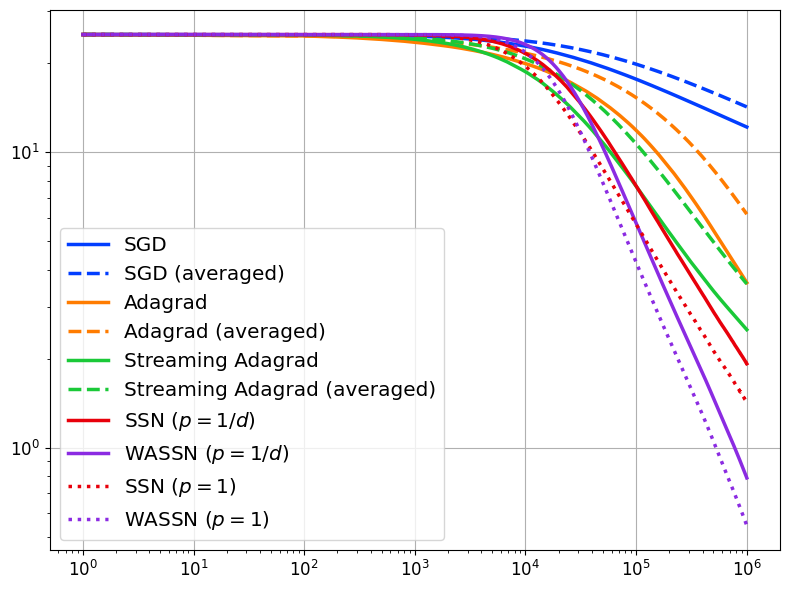}
     \caption{Logistic Regression: Mean-squared error of the distance to the optimum $\theta^{*}$, plotted against the sample size of $1.000.000$, for various initializations. The initial points $\theta_{0}$ are generated as $\theta_{0}=\theta^{*}(1+rU)$, where $U$ is a uniform random variable on the unit sphere of $\R^{d}$, and $r$ takes values of $1$ (left) or $5$ (right). Each curve reports $\lVert\theta_{t}-\theta^{*}\rVert$ averaged over $50$ different epochs, with a different initial point drawn for each sample.}
    \label{fig:logistic}
\end{figure*}

In \cref{fig:logistic}, we display the evolution of the quadratic mean error of different estimates, for two different initializations. Across all initial configurations, the stochastic Newton's method consistently emerges as a strong competitor. Conversely, the effectiveness of the Adagrad algorithm diminishes as the initial starting point moves further away from the solution. In scenarios involving less-than-ideal initializations, as depicted on the right side of the figure, the best performances are achieved by the averaged Newton's method. This exceptional asymptotic behavior is enabled by the incorporation of weighted estimates, assigning greater significance to the most recent ones, distinguishing it from the "standard" averaging Newton's method, as elaborated in \citep{bercu2020efficient}.

\section*{Conclusions and Future Work} \label{sec::con}
In this work, we addressed the unique challenges posed by streaming data in the context of stochastic optimization. The continuous influx of large, high-dimensional data necessitates adaptive approaches that can effectively handle ill-conditioned problems while maintaining computational efficiency. Our contributions lie in the development of adaptive stochastic optimization methods, particularly an inversion-free adaptive Newton's method with a computational complexity matching that of first-order methods, $\mathcal{O}(dN_{t})$, where $d$ represents the number of dimensions/features, and $N_{t}$ denotes the quantity of data up to time $t$.

\medskip

Theoretical analyses have confirmed the asymptotic efficiency of our proposed methods. By dynamically adjusting learning rates per-dimension and incorporating historical gradient or Hessian information, our methods exhibit adaptability and efficiency in navigating through the complexities of ill-conditioned problems. Notably, the introduction of a weighted averaged version enhances the adaptability and robustness of our methods, particularly in scenarios involving complex covariance structures and challenging initializations.

\medskip

One significant contribution is the inversion-free adaptive Newton's method in \cref{sec::app::newton:nd:wa}, which strikes a balance between addressing ill-conditioned problems and meeting the computational demands of streaming data. This innovation allows us to harness the advantages of second-order information while aligning with the computational complexity of first-order methods. Empirical evidence demonstrates the effectiveness of our adaptive methods, showcasing superior performance, especially in challenging scenarios.

\medskip

In conclusion, our adaptive stochastic optimization methods offer a versatile solution for streaming data settings, providing an efficient and adaptive framework for handling ill-conditioned problems. The inversion-free adaptive Newton's method, in particular, stands out as a computationally efficient alternative that bridges the gap between first-order and second-order methods. As we look ahead, further exploration of real-world applications, theoretical advancements, and extensions to non-convex settings will be key directions for future research in this evolving field.

\medskip

\textbf{Future work.} Looking ahead, there are several promising directions for future research: (a) Non-convex analysis: Extending our methodologies to non-convex optimization problems is a crucial next step. Analyzing the behavior and convergence properties of our adaptive methods in non-convex scenarios will contribute to a more comprehensive understanding of their applicability across diverse optimization landscapes. (b) Time-dependent observations: The streaming context often involves time-dependent observations, and our current work assumes independence among the data points. Investigating extensions of our methods to handle dependent observations will be essential for real-world applications where temporal or spatial dependencies are prevalent. Recently, \citet{godichon2023learning} showed that increasing mini-batches can break both short- and long-term dependence structures. These future research directions aim to refine the versatility and robustness of our adaptive stochastic optimization methods, ensuring their effectiveness across a broader spectrum of optimization challenges.


\subsection*{Aknowledgements}
{N. Werge acknowledges the support of the Novo Nordisk Foundation (NNF) through grant NNF21OC0070621.}

\appendix
\crefalias{section}{appendix}
\setcounter{equation}{0}
\renewcommand{\theequation}{\thesection.\arabic{equation}}
\setcounter{thm}{0}
\renewcommand{\thethm}{\thesection.\arabic{thm}}
\setcounter{lem}{0}
\renewcommand{\thelem}{\thesection.\arabic{lem}}
\setcounter{prop}{0}
\renewcommand{\theprop}{\thesection.\arabic{prop}}
\setcounter{cor}{0}
\renewcommand{\thecor}{\thesection.\arabic{cor}}
\setcounter{ass}{0}
\renewcommand{\theass}{\thesection.\arabic{ass}}
\setcounter{ex}{0}
\renewcommand{\theex}{\thesection.\arabic{ex}}
\setcounter{rem}{0}
\renewcommand{\therem}{\thesection.\arabic{rem}}

\section*{Appendix}
\Cref{sec::appendix} contains the statements for increasing mini-batches, while \cref{sec::appendix:proofs} contains the mathematical proofs of the main results.

\section{Statements for Increasing Mini-Batches} \label{sec::appendix}
In this appendix, we investigate our adaptive methods, in the case, where the mini-batches are increasing, and give the translation of the different theorems in this case. Following, \citet{godichon2023learning,godichon2023non}, we consider mini-batch sizes of the form $n_{t}= \lfloor C_{\rho}t^{\rho} \rfloor$ with $C_{\rho}\in\N$. In this case, we take $\gamma_{t} = C_{\gamma}n_{t}^{\beta}t^{-\gamma}$, which roughly means that $\gamma_{t} \sim C_{\gamma}C_{\rho}^{\beta}t^{-\gamma + \beta \rho}$. Adding the term $n_{t}^{\beta}$ to the learning rate enables us to give more weights on (presumably) more precise gradient steps, as they are estimated with larger mini-batches $n_{t}$. We suppose that $\gamma - \beta \rho \in (1/2,1)$ and $\gamma > \frac{\rho (2\beta -1) +1 }{2}$. 

\subsection{Adaptive Stochastic Optimization Methods} \label{sec::appendix::asom}
When considering increasing mini-batches, our adaptive stochastic optimization methods are defined as in \cref{eq::asom}:
\begin{equation} \label{eq:asom:increasing}
\theta_{t+1} = \theta_{t} - \gamma_{t+1} A_{t} \nabla_{\theta} f(\theta_{t};\xi_{t+1}), \; \theta_{0}\in\R^{d},
\end{equation}
where $\nabla_{\theta}f(\theta_{t};\xi_{t+1}) = n_{t+1}^{-1}\sum_{i=1}^{n_{t+1}} \nabla_{\theta} f( \theta_{t};\xi_{t+1,i})$. In this case, the conditions in \cref{cond::step} should be rewritten to
\begin{align}\label{cond::step::increasing}
\sum_{t \geq 1}\gamma_{t} \lambda_{\min}(A_{t-1}) = + \infty \; \as, \quad \sum_{t \geq 1} \frac{\gamma_{t}^{2}}{n_{t}} \lambda_{\max} ( A_{t-1})^{2} < + \infty \; \as, \quad \frac{\lambda_{\max}(A_{t})^{2}\gamma_{t+1}}{\lambda_{\min}(A_{t})} \xrightarrow[t\to + \infty]{\as} 0.
\end{align}

Under these conditions in \cref{cond::step::increasing}, we can show the strongly consistency of the estimates derived from \cref{eq:asom:increasing}. This, in the increasing mini-batch case, \cref{theo::ps} can be rewritten as follows:
\begin{thm}\label{theo::ps::increasing}
Suppose \cref{ass::1,ass::2,ass::3} hold, along with the conditions in \cref{cond::step::increasing}. Then $\theta_{t}$ converges almost surely to $\theta$.
\end{thm}

Similarly, we have the rate of convergence for \cref{eq:asom:increasing} as in \cref{theo::rate}:
\begin{thm}\label{theo::rate::increasing}
Suppose \cref{ass::1,ass::2,ass::3,ass::4} hold, along with the conditions in \cref{cond::step::increasing}. In addition, assume there exists positive constants $C_{\eta}$ and $\eta > \frac{1}{\gamma - \beta \rho}-1$ such that for all $\theta\in\R^{d}$,
\begin{equation} \label{momenteta::increasing}
\E\left[ \lVert \nabla_{\theta}f(\theta;\xi) \rVert^{2+2\eta}\right] \leq C_{\eta} ( 1+ F(\theta_{t}) -F(\theta^{*}) )^{1+\eta}.
\end{equation}
Then, 
\begin{equation*}
\lVert \theta_{t} - \theta^{*} \rVert^{2} = \mathcal{O} \left( \ln(N_{t})  N_{t}^{\frac{-\rho - \gamma + \beta \rho }{1+\rho}}  \right) \; \as
\end{equation*}
\end{thm}
Note that the rate of convergence in \cref{theo::rate::increasing} reproduce the results of the constant mini-batch case in \cref{theo::rate} when $n_{t}=n=C_{\rho}$, $\beta=0$, and $\rho=0$.

\subsection{The Weighted Averaged Version} \label{sec::appendix::asom::wa}
As in \cref{sec::appendix::asom}, we consider the weighted Polyak-Ruppert averaged version of our adaptive stochastic optimization methods for increasing mini-batches; the constant mini-batch case is in \cref{sec::asom::wa}. These weighted estimates are defined for $w\geq0$ as follows:
\begin{equation*}
\theta_{t,w}= \frac{1}{\sum_{i=0}^{t-1}n_{i+1}\ln(i+1)^{w}} \sum_{i=0}^{t-1}n_{i+1}\ln(i+1)^{w} \theta_{i},
\end{equation*}
which can be written recursively as
\begin{equation*}
\theta_{t+1,w} = \theta_{t,w} + \frac{n_{t+1} \ln(t+1)^{w}}{\sum_{i=0}^{t}n_{i+1}\ln(i+1)^{w}} (\theta_{t}-\theta_{t,w}).
\end{equation*}

Likewise to \cref{theotlc}, we have the rate of convergence and the optimal asymptotic normality of these weighted estimates:
\begin{thm}\label{theotlc::increasing}
Suppose \cref{ass::1,ass::2,ass::3,ass::4,ass::5} hold, along with inequality \cref{momenteta::increasing}. In addition, assume there exists a positive constant $\nu$ such that
\begin{equation*}
\lVert A_{t} - A \rVert_{\op} = \mathcal{O} \left( \frac{1}{t^{\nu}} \right) \; \as
\end{equation*}
Then,
\begin{equation*}
 \lVert \theta_{t,w} - \theta^{*} \rVert^{2} = \left\lbrace \begin{array}{lr}  \mathcal{O} \left( \frac{ \ln(N_{t})^{\frac{1}{2}+ \mathbb{1}_{\left\{\nu + \frac{\rho(1-\beta) + \gamma}{2}=1\right\}}}}{N_{t}^{\frac{  2\nu + \rho(1-\beta) + \gamma}{1+\rho}}} \right) \; \as,
& \text{if } \;   2\nu + \rho(1-\beta) + \gamma \leq 1+\rho, \\
\mathcal{O} \left( \frac{\ln(N_{t})}{N_{t}} \right) \;  \as, & \text{if } \; 2\nu + \rho(1-\beta) + \gamma > 1+\rho.
\end{array} \right.
\end{equation*}
Moreover, if $2\nu + \rho(1-\beta) + \gamma > 1+\rho$, then
\begin{equation*}
\sqrt{N_{t}} ( \theta_{t,w} - \theta^{*} ) \xrightarrow[t\to +\infty]{\mathcal{L}} \mathcal{N} ( 0 ,\nabla_{\theta}^{2}F(\theta^{*})^{-1}\Sigma \nabla_{\theta}^{2}F(\theta^{*})^{-1}) .
\end{equation*}
\end{thm}

Similarly to \cref{theowithoutrate}, we can establish the asymptotic efficiency without relying on a (weak) rate of convergence of $A_{t}$:
\begin{thm} \label{theowithoutrate::increasing}
Suppose \cref{ass::1,ass::2,ass::3,ass::4,ass::5} hold, along with inequality \cref{momenteta::increasing}. In addition, assume there exists a positive constant $v'> 1/2$ such that
\begin{equation}\label{equalitystrange::increasing}
\frac{1}{\sum_{i=0}^{t-1}n_{i+1} \ln(i+1)^{w}} \sum_{i=0}^{t-1} n_{i+1} \ln(i+1)^{w+1/2+\delta} \lVert A_{i+1}^{-1} - A_{i}^{-1} \rVert_{\op} (i+1)^{\frac{\gamma-\rho(\beta +1)}{2}} = \mathcal{O} \left( \frac{1}{t^{(1+\rho )v'}} \right) \; \as,
\end{equation}
for some $\delta>0$. Then,
\begin{equation*}
\lVert \theta_{t,w} - \theta^{*} \rVert^{2} = \mathcal{O} \left( \frac{\ln(N_{t})}{N_{t}} \right) \; \as, \quad \text{and} \quad \sqrt{N_{t}} (\theta_{t,w} - \theta^{*} ) \xrightarrow[t\to + \infty]{\mathcal{L}} \mathcal{N} ( 0 , \nabla_{\theta}^{2}F(\theta^{*})^{-1}\Sigma \nabla_{\theta}^{2}F(\theta^{*})^{-1}) .
\end{equation*}
\end{thm}

\subsection{Applications to Newton's Methods}
As in \cref{sec::app}, we apply our adaptive stochastic optimization methodology to (stochastic) Newton's methods (but with increasing mini-batches). In particular, we consider the Newton's methods with the possibly $\mathcal{O}(dN_{t})$ operations, analogues to \cref{sec::app::newton:nd,sec::app::newton:nd:wa}.

\subsubsection{Streaming Stochastic Newton's Methods with possibly $\mathcal{O}(dN_{t})$ operations}
Expanding the mini-batch scenario from \cref{sec::app::newton:nd} leads to the formulation of the streaming variant of stochastic Newton's method, as defined by:
\begin{equation} \label{def::newton::increasing}
\theta_{t+1} = \theta_{t} - \frac{1}{N_{t+1}} \bar{H}_{t,w'}^{-1} \sum_{i=1}^{n_{t+1}} \nabla_{\theta}f (\theta_{t};\xi_{t+1,i}),
\end{equation}
where $\bar{H}_{t,w'}=N_{t,Z}^{-1}H_{t,w'}$ with 
\begin{equation*}
H_{t,w'}  = H_{0,w'}  +  \sum_{t'=1}^{t}\ln (t'+1)^{w'} \sum_{i=1}^{n_{t}'}Z_{t',i} \left(  \iota_{t',i} \tilde{e}_{t',i}\tilde{e}_{t',i}^{T} +   \alpha_{t',i} \Phi_{t',i}\Phi_{t',i}^{\top} \right),
\end{equation*}
with $N_{Z,t} = 1+ \sum_{t'=1}^{t}\ln (t'+1)^{w'} \sum_{i=1}^{n_{t'}} Z_{t',i}$.
In addition, let $N_{Z,t,i} = (1+ \sum_{t'=1}^{t-1}\sum_{j=1}^{n_{t'}}Z_{t',j} + \sum_{j=1}^{i} Z_{t,j})$,  $\iota_{t',i} = c_{\iota} N_{Z,t',i} ^{-\iota}$,  $\iota \in (0 \frac{1-\rho}{2(1+\rho)})$, $e_{t',i} $ is the ($N_{t'-1} +i$ modulo $d$ +1)-th componant of the canonical basis.

As in \cref{theo::newton::nd}, we can establish the rate of convergence and the asymptotic normality of \cref{def::newton::increasing}:
\begin{thm}\label{theo::newton::nd::increasing}
Suppose \cref{ass::1,ass::2,ass::3,ass::5} hold, along with the conditions in \cref{momenteta::increasing,hessegale}. In addition, assume there are positive constants $C_{\eta'}$ and $\eta'>1$ such that for all $\theta\in\R^{d}$,
\begin{equation*}
\E \left[ \lVert \alpha(\theta;\xi)\Phi(\theta;\xi)\Phi(\theta;\xi)^{\top} \rVert^{\eta'} \right] \leq C_{\eta'}^{\eta'}.
\end{equation*}
Then, 
\begin{equation*}
\lVert \theta_{t} - \theta^{*} \rVert^{2} = \mathcal{O} \left( \frac{\ln(N_{t})}{N_{t}} \right) \; \as
\end{equation*}
Moreover, suppose that the Hessian of $F$ is locally Lipschitz on a neighborhood around $\theta^{*}$ and that $\eta'\geq 2$. Then,
\begin{equation*}
\sqrt{N_{t}} \left( \theta_{t} - \theta^{*} \right) \xrightarrow[t \to + \infty]{\mathcal{L}}\mathcal{N}\left( 0 , \nabla_{\theta}^{2}F(\theta^{*})^{-1}\Sigma \nabla_{\theta}^{2}F(\theta^{*})^{-1} \right).
\end{equation*}
\end{thm}

\subsubsection{Weighted Averaged Version of Streaming Stochastic Newton's Methods with possibly $\mathcal{O}(dN_{t})$ operations}
The weighted averaged version outlined in \cref{sec::app::newton:nd:wa} can similarly be adapted to the increasing mini-batch case. The weighted averaged streaming stochastic Newton's method is defined as
\begin{align}
& \notag \theta_{t+1} = \theta_{t} -\gamma_{t+1} \bar{S}_{t,w'}^{-1} \frac{1}{n_{t+1}} \sum_{i=1}^{n_{t+1}} \nabla_{\theta}f(\theta_{t}; \xi_{t+1,i}), \\
& \label{def::WASN::increasing} \theta_{t+1,w} = \theta_{t,w} + \frac{n_{t+1}\ln(t+1)^{w}}{\sum_{i=0}^{t}n_{i+1}\ln(i+1)^{w}}( \theta_{t} - \theta_{t,w}),
\end{align}
where $\gamma_{t} = C_{\gamma}n_{t}^{\beta}t^{-\gamma}$ and $\bar{S}_{t,w'}=N_{t,Z}^{-1}S_{t,w'}$ with
\begin{equation*}
S_{t,w'} = S_{0,w'} + \sum_{t'=1}^{t}\ln(t'+1)^{w'} \sum_{i=1}^{n_{t'}} Z_{t',i}\left(  \iota_{t',i}  e_{t',i}e_{t',i}^{T} +  \alpha_{t',i} \Phi_{t',i}\Phi_{t',i}^{\top}  \right),
\end{equation*}
with $S_{0}$ symmetric and positive, $\iota_{t',i} = C_{\iota} N_{Z,t',i} ^{-\iota}$ with $
\iota \in \left(0,  \frac{\min \left\lbrace \gamma - \rho \beta , 2 \gamma -2 \rho \beta -1 + \rho \right\rbrace}{2(1+\rho)}\right)$,
which is possible since $\gamma -  \rho \beta \in (1/2,1)$.

Like in \cref{theo::newton::wasn}, we have the following asymptotic optimality:
\begin{thm} \label{theo::newton::wasn::inceasing}
Suppose \cref{ass::1,ass::2,ass::3,ass::5} hold, along the conditions in \cref{momenteta::increasing,hessegale}. In addition, assume there exists positive constants $C_{\eta'}$ and $\eta'>1$ such that for all $\theta\in\R^{d}$,
\begin{equation*}
\E\left[ \lVert \alpha(\theta;\xi)\Phi(\theta;\xi) \Phi(\theta;\xi)^{\top} \rVert^{\eta'} \right] \leq C_{\eta'}^{\eta'}.
\end{equation*}
Then, 
\begin{equation*}
\lVert \theta_{t} - \theta^{*} \rVert^{2} = \mathcal{O}\left( \frac{\ln(N_{t})}{N_{t}^{\frac{\gamma + \rho (1-\beta )}{1+\rho}}} \right) \; \as \qquad \text{and} \qquad \lVert \theta_{t,w} - \theta^{*} \rVert^{2} = \mathcal{O} \left( \frac{\ln(N_{t})}{N_{t}} \right).
\end{equation*}
In addition,
\begin{equation*}
\sqrt{N_{t}} \left( \theta_{t,w} - \theta^{*} \right) \xrightarrow[t \to + \infty]{\mathcal{L}}\mathcal{N}\left( 0 , \nabla_{\theta}^{2}F(\theta^{*})^{-1}\Sigma \nabla_{\theta}^{2}F(\theta^{*})^{-1} \right) .
\end{equation*}
\end{thm}

\subsection{Streaming Adagrad and its Weighted Averaged Version} \label{sec::app::adagrad}
In this section, we apply our adaptive stochastic optimization methodology to Adagrad \citep{duchi2011adaptive}. Our adaptation results in a streaming version of Adagrad, specifically tailored for efficient handling of evolving data streams. Additionally, we introduce the weighted averaged version of streaming Adagrad, enhancing adaptability and accelerating convergence.

\subsubsection{Streaming Adagrad with constant mini-batches}  The recursive definitions for streaming Adagrad and its weighted averaged version are as follows:
\begin{align} \label{eq::adagrad}
\theta_{t+1} =& \theta_{t} - \gamma_{t+1} G_{t} \nabla_{\theta} f(\theta_{t};\xi_{t+1}), \quad \theta_{0}\in\R^{d}, \\
\theta_{t+1,w} =& \theta_{t,w} + \frac{\ln(t+1)^{w}}{\sum_{i=0}^{t}\ln(i+1)^{w}}(\theta_{t} - \theta_{t,w} ),
\end{align}
where $\nabla_{\theta}f(\theta_{t};\xi_{t+1})=n^{-1}\sum_{i=1}^{n} \nabla_{\theta}f(\theta_{t};\xi_{t+1,i})$ and $G_{t}$ is a diagonal matrix with $k$-th element $G_{t}^{(k)}$ for $k=1,\dots,d$, given as
\begin{equation*}
G_{t}^{(k)} = \left(\frac{1}{N_{t}}\left( G_{0}^{(k)} + \sum_{i=1}^{t}\sum_{j=1}^{n}  \frac{\partial}{\partial k}f(\theta_{t-1};\xi_{i,j})^{2}   \right)\right)^{-1/2},
\end{equation*}
with $\nabla_{\theta^{(k)}}$ denoting the partial derivative with respect to $k$-th element of $\theta$, i.e., $\theta^{(k)}$.

To mitigate the potential divergence of the eigenvalues of $G_{t}$, we employ a technique introduced by \citet{GBT2023}, resulting in a mild modification of the standard random matrix $G_{t}$. The modification is expressed as:
\begin{equation*}
G_{t}^{(k)} = \max \left\lbrace C_{\beta''}t^{\beta''} , \min \left\lbrace C_{\beta'}t^{\beta'} , \left(\frac{1}{N_{t}}\left( G_{0}^{(k)} + \sum_{i=1}^{t}\sum_{j=1}^{n} \frac{\partial}{\partial k}f(\theta_{t-1};\xi_{i,j})^{2}   \right)\right)^{-1/2} \right\rbrace \right\rbrace,
\end{equation*}
with $C_{\beta'},C_{\beta''}>0$. In this formulation, the addition of the $\min$-term in $G_{t}$ aids in controlling the potential divergence of its largest eigenvalue, while the $\max$-term ensures a lower bound for the smallest eigenvalue. Precisely, selecting $\gamma\in(1/2,1)$, $\beta'\in(0,\gamma-1/2)$, and $\beta''\in(\gamma-1,0)$ satisfies $2\beta'-\gamma-\beta'' < 0$, which ensures the conditions in \cref{cond::step} are satisfied.

With these modifications in place, we can now establish the rate of convergence and asymptotic normality.
\begin{thm} \label{theo::adagrad}
Suppose \cref{ass::1,ass::2,ass::3,ass::5} hold, along with inequality \cref{eq::momenteta}.  {In addition, assume that   the variance $\V[\frac{\partial}{\partial k}f(\theta^{*};\xi)]>0$ for $k=1,\dots,d$.} Then,  
\begin{equation*}
\lVert \theta_{t} - \theta^{*} \rVert^{2} = \mathcal{O}\left( \frac{\ln(N_{t})}{N_{t}^{\gamma}} \right) \; \as, \qquad \lVert \theta_{t,w} - \theta^{*} \rVert^{2}  = \mathcal{O} \left(\frac{\ln(N_{t})}{N_{t}} \right) \; \as,
\end{equation*}
and
\begin{equation*}
\sqrt{N_{t}}(\theta_{t,w} - \theta^{*} ) \xrightarrow[t\to+ \infty]{\mathcal{L}} \mathcal{N}(0,\nabla_{\theta}^{2}F(\theta^{*})^{-1}\Sigma\nabla_{\theta}^{2}F(\theta^{*})^{-1}).
\end{equation*}
\end{thm}

\subsubsection{Streaming Adagrad with increasing mini-batches}
For the increasing mini-batch case, the streaming Adagrad variant and its weighted averaged version is defined recursively by
\begin{align*}
\theta_{t+1}  & = \theta_{t}  - \gamma_{t+1} G_{t} \nabla_{\theta} f(\theta_{t}; \xi_{t+1}), \theta_{0}\in\R^{d}, \\
\theta_{t+1,w} & = \theta_{t,w}  + \frac{n_{t+1} \ln (t+1)^{w}}{\sum_{i=0}^{t} n_{i+1}\ln(i+1)^{w}} (\theta_{t}-\theta_{t,w}),  
\end{align*}
where $\nabla_{\theta}f(\theta_{t};\xi_{t+1})=n_{t+1}^{-1}\sum_{i=1}^{n_{t+1}} \nabla_{\theta}f(\theta_{t};\xi_{t+1,i})$ and $G_{t}$ is a diagonal matrix with, denoting by $G_{t}^{(k)}$ the $k$-th element of the diagonal of $G_{t}$,
\begin{equation*}
G_{t}^{(k)} = \max \left\lbrace C_{\beta''}t^{\beta''} , \min \left\lbrace C_{\beta'}t^{\beta'} , \left( \frac{1}{N_{t}}\left( G_{0}^{(k)}  + \sum_{i=1}^{t}\sum_{j=1}^{n_{t}} \left( \frac{\partial}{\partial k}f(\theta_{t-1};\xi_{i,j})\right)^{2}   \right)\right)^{-1/2} \right\rbrace \right\rbrace .
\end{equation*}
Remark that the add of the minimum in the expression of $G_{t}$ enables to control the possible divergence of the larges eigenvalue of $G_{t}$ while the max term enables to lower bound the smallest eigenvalue. More precisely, taking $\gamma -\beta \rho \in (1/2,1)$, $\beta ' \in (0 , \gamma+ \rho(\frac{1}{2} - \beta ) -1/2)$ and $\beta '' \in (\gamma - \beta \rho -1 , 0)$  satisfying $2\beta ' - \gamma + \beta \rho - \beta ''  < 0$ enables to verify the  conditions in \cref{cond::step::increasing}. To simplify it, one can take $\beta ' < \gamma - \beta \rho - \nicefrac{1}{2}$.   
Then, \cref{theo::adagrad} can be written as follows:
\begin{thm} \label{theo::adagrad::increasing}
Suppose \cref{ass::1,ass::2,ass::3,ass::5} hold, along with the conditions in \cref{momenteta::increasing}. In addition, assume that the variance $\mathbb{V} \left[ \frac{\partial}{\partial k} f(\theta^{*};\xi) \right]>0$ for $k=1,\dots,d$. Then,  
\begin{equation*}
\lVert \theta_{t} - \theta^{*} \rVert^{2} = \mathcal{O} \left( \frac{\ln(N_{t})}{N_{t}^{\frac{\gamma +  \rho (1-\beta)}{1+\rho}}} \right) \; \as, \qquad \lVert \theta_{t,w} - \theta^{*} \rVert^{2} = \mathcal{O} \left( \frac{\ln(N_{t})}{N_{t}} \right) \; \as,
\end{equation*}
and
\begin{equation*}
\sqrt{N_{t}}( \theta_{t,w} - \theta^{*}) \xrightarrow[t\to + \infty]{\mathcal{L}} \mathcal{N} \left( 0 , \nabla_{\theta}^{2}F(\theta^{*})^{-1}\Sigma \nabla_{\theta}^{2}F(\theta^{*})^{-1} \right).
\end{equation*}
\end{thm}

\section{Proofs} \label{sec::appendix:proofs}
The proof are solely presented for the increasing mini-batch case outlined in \cref{sec::appendix}, as the constant mini-batch case corresponds to $n_{t}=n=C_{\rho}$, $\beta=0$, and $\rho = 0$.

For the sake of simplicity, in all the sequel, since $n_{t} \sim C_{\rho}t^{\rho}$, we will make the abuse that $n_{t} = C_{\rho} {t}^{\rho}$. To lighten the notation, we let $H$ denote $\nabla_{\theta}^{2}F(\theta^{*})$. In addition, let $f_{t+1}'$ and $f_{t+1,i}'$ denote $\nabla_{\theta}f(\theta_{t};\xi_{t+1})$ and $\nabla_{\theta}f(\theta_{t};\xi_{t+1,i})$, respectively.

\subsection{Proof of \cref{theo::ps,theo::ps::increasing}}
Let $V_{t}$ denote $F(\theta_{t}) - F(\theta^{*})$.
Observe that with the help of a Taylor's expansion of the objective function $F$ and since the Hessian is uniformly bounded (\cref{ass::2}), then one has 
\begin{align*}
V_{t+1} & \leq  V_{t} + \nabla_{\theta} F(\theta_{t})^{\top}( \theta_{t+1} - \theta_{t})   + \frac{L_{\nabla F}}{2} \lVert \theta_{t+1} - \theta_{t} \rVert^{2}  \\
& \leq V_{t} - \gamma_{t+1}\nabla_{\theta} F(\theta_{t})^{\top} A_{t}  f_{t+1}'  + \frac{L_{\nabla F}}{2} \gamma_{t+1}^{2} \lambda_{\max} (A_{t})^{2} \lVert f_{t+1}' \rVert^{2}.
\end{align*}
Before taking the conditional expectation, recall from \citet{godichon2023non} that
\begin{equation*}
\E [  \lVert f_{t+1}' \rVert^{2} \vert \mathcal{F}_{t} ] = \frac{1}{n_{t+1}} \E[ \lVert f_{t+1,i}' \rVert^{2} \vert \mathcal{F}_{t} ] + \lVert \nabla_{\nabla} F(\theta_{t}) \rVert^{2} \leq   \frac{1}{n_{t+1}} C(1+F(\theta_{t})-F(\theta^{*})) + \lVert \nabla_{\theta} F(\theta_{t}) \rVert^{2}.
\end{equation*}
Thus, we obtain that
\begin{align*}
\E [ V_{t+1} \vert \mathcal{F}_{t}] \leq & V_{t} - \gamma_{t+1}\nabla_{\theta} F (\theta_{t})^{\top} A_{t} \nabla_{\theta} F(\theta_{t}) \\ & + \frac{L_{\nabla F}}{2} \gamma_{t+1}^{2} \lambda_{\max}(A_{t})^{2} \left( \frac{C}{n_{t+1}} \left( 1+ F(\theta_{t}) -  F(\theta^{*}) \right) + \lVert \nabla_{\theta} F(\theta_{t}) \rVert^{2} \right) \\
\leq & \left( 1+ \frac{L_{\nabla F}C}{2} \frac{\gamma_{t+1}^{2}\lambda_{\max}(A_{t})^{2}}{n_{t+1}} \right) V_{t} - \gamma_{t+1}\lVert \nabla_{\nabla} F(\theta_{t}) \rVert^{2}  \left( \lambda_{\min}(A_{t})  - \frac{L_{\nabla F}}{2} \gamma_{t+1}\lambda_{\max}(A_{t})^{2} \right)  \\
& + \frac{L_{\nabla F}C}{2}\frac{\gamma_{t+1}^{2} \lambda_{\max} (A_{t})^{2}}{n_{t+1}}.
\end{align*}
Observe that as $\frac{\gamma_{t+1}\lambda_{\max}(A_{t})^{2}}{\lambda_{\min}(A_{t})}$ converges almost surely to zero for any constant $C \in (0,1)$ due to the conditions in \cref{cond::step::increasing}. Then, $\mathbb{1}_{ \left\lbrace \frac{L_{\nabla F} }{2} \gamma_{t+1}\lambda_{\max} \left( A_{t} \right)^{2} \geq  C \lambda_{\min} \left( A_{t} \right) \right\rbrace} $ converges almost surely to zero as well. Thus, we have that
\begin{align*}
\E [V_{t+1}\vert\mathcal{F}_{t}] \leq & \left( 1+ \frac{L_{\nabla F}C}{2}\frac{\gamma_{t+1}^{2}\lambda_{\max}(A_{t})^{2}}{n_{t+1}} \right) V_{t} - (1-C)\gamma_{t+1} \lambda_{\min}(A_{t}) \lVert \nabla_{\theta} F(\theta_{t}) \rVert^{2}  \\ & + \frac{L_{\nabla F}C}{2}\frac{\gamma_{t+1}^{2} \lambda_{\max}(A_{t})^{2}}{n_{t+1}} \\ & +  \frac{L_{\nabla F}C}{2}\gamma_{t+1}^{2}\lambda_{\max}(A_{t})^{2} \lVert\nabla_{\theta} F(\theta_{t}) \rVert^{2} \mathbb{1}_{\left\lbrace \frac{L_{\nabla F} }{2} \gamma_{t+1}\lambda_{\max}(A_{t})^{2} \geq C \lambda_{\min}(A_{t}) \right\rbrace}.
\end{align*}
Next, since $\mathbb{1}_{\left\lbrace \frac{L_{\nabla F} }{2}\gamma_{t+1}\lambda_{\max}(A_{t})^{2} \geq  C\lambda_{\min}(A_{t}) \right\rbrace}$ converges almost surely to zero and by the conditions in \cref{cond::step::increasing};
\begin{equation*}
\sum_{t \geq 0} \frac{\gamma_{t+1}^{2}\lambda_{\max} \left( A_{t} \right)^{2}}{n_{t+1}}  < + \infty \; \as, 
\end{equation*}
and 
\begin{equation*}
\sum_{t\geq 0}  \gamma_{t+1}^{2}\lambda_{\max} \left( A_{t} \right)^{2} \mathbb{1}_{\left\lbrace \frac{L_{\nabla F}C}{2}\gamma_{t+1} \left\| \nabla F \left( \theta_{t} \right) \right\|^{2}\lambda_{\max} \left( A_{t} \right)^{2} \geq c \lambda_{\min} \left( A_{t} \right) \right\rbrace} < + \infty \; \as,
\end{equation*}
then, applying Robbins-Siegmund's theorem gives that $V_{t}$ converges almost surely to a finite random variable and
\begin{equation*}
\sum_{t \geq 0} \gamma_{t+1}  \lambda_{\min}(A_{t}) \lVert \nabla_{\theta} F(\theta_{t}) \rVert^{2} < + \infty \; \as,
\end{equation*}
meaning, that $\liminf_{t}\lVert \nabla_{\theta} F(\theta_{t}) \rVert^{2}=0$ a.s., such that $\liminf_{t} V_{t} = 0$ a.s., i.e., $V_{t}$ converges almost surely to zero, which concludes the proof.

\subsection{Proof of \cref{theo::rate,theo::rate::increasing}}
Following the reasoning of \citet[page 11]{antonakopoulos2022adagrad}, $AH$ and $A^{1/2}HA^{1/2}$ have the same eigenvalues. Indeed, for any $\lambda \in \mathbb{R}$,
\begin{align*}
\det \left( A^{1/2}HA^{1/2} - \lambda I_{d} \right) &  = \det \left( A^{-1/2} \left( AHA^{1/2} - \lambda A^{1/2} \right) \right) \\
& = \det \left( A^{-1/2} \left( AH - \lambda I_{d} \right) A^{1/2} \right) \\
&  = \det \left( AH - \lambda I_{d} \right) .
\end{align*}
Then, there exists matrix $Q$ and a positive diagonal matrix $D$, such that $AH = Q^{-1}D Q$. Thus,  
\begin{equation*}
Q(\theta_{t+1} - \theta^{*}) = Q( \theta_{t} - \theta^{*}) - \gamma_{t+1}QA_{t} f_{t+1}' = Q( \theta_{t} - \theta^{*}) - \gamma_{t+1}QA_{t} \nabla_{\theta} F( \theta_{t}) + \gamma_{t+1}QA_{t} \xi_{t+1},
\end{equation*}
where $\xi_{t+1} = \nabla_{\theta} F( \theta_{t}) - f_{t+1}'$. By linearizing the gradient one has
\begin{align*}
Q( \theta_{t+1} - \theta^{*}) = Q( \theta_{t} - \theta^{*}) - \gamma_{t+1} QA_{t} H( \theta_{t} - \theta^{*} ) + \gamma_{t+1}QA_{t} \xi_{t+1} - \gamma_{t+1}QA_{t} \delta_{t},
\end{align*}
where $\delta_{t} = \nabla_{\theta} F(\theta_{t}) - H( \theta_{t} - \theta^{*})$ is the remainder term of the Taylor's expansion of the gradient. Next, we have
\begin{align}
 \notag    Q \left( \theta_{t+1} - \theta^{*} \right) = & Q \left( \theta_{t} - \theta^{*} \right)- \gamma_{t+1} QAH \left( \theta_{t} - \theta^{*} \right) - \gamma_{t+1} Q\left( A_{t} - A \right) H \left( \theta_{t} - \theta^{*} \right) \\ \notag & + \gamma_{t+1}QA_{t} \xi_{t+1} - \gamma_{t+1}QA_{t} \delta_{t} \\
 \notag = & \left( I_{d} - \gamma_{t+1}D \right) Q \left( \theta_{t} - \theta^{*} \right)  - \gamma_{t+1}Q\left( A_{t} - A \right) H \left( \theta_{t} - \theta^{*} \right) \\  & + \gamma_{t+1}QA_{t} \xi_{t+1} - \gamma_{t+1}QA_{t} \delta_{t}. \label{decdelta} 
\end{align}
Observe that in the case where $A=H^{-1}$, i.e., in the stochastic Newton's method, one has $D= Q= I_{d}$. With the help of induction, one has by \cref{decdelta} that
\begin{align}
\notag   Q \left( \theta_{T} - \theta^{*} \right)&  = \overbrace{ \beta_{T,0} Q \left( \theta_{0} - \theta^{*} \right) }^{R_{1,T}:=} \overbrace{ - \sum_{t=0}^{T-1} \beta_{T,t+1} \gamma_{t+1}Q \left( A_{t} - A \right) H \left( \theta_{t} - \theta^{*} \right) -  \sum_{t=0}^{T-1} \beta_{T,t+1} \gamma_{t+1}QA_{t} \delta_{t}}^{R_{2,T}:=} \\
\label{decbeta} &  + \underbrace{ \sum_{t=0}^{T-1} \beta_{T,t+1} \gamma_{t+1}QA_{t} \xi_{t+1}}_{M_{T}:=}  ,
\end{align}
where $\beta_{T,t} = \prod_{j=t+1}^{T}\left( I_{d} - \gamma_{j} D \right) $ and $\beta_{T,T} = I_{d}$.
The rest of the proof consists in giving the rate of convergence of each term on the right-hand side of decomposition \cref{decbeta} for both cases, i.e., for the constant mini-batches and increasing mini-batches.

\paragraph*{Rate of convergence for $R_{1,T}$.} Since $D$ is a positive diagonal matrix, and since $\gamma_{t}$ is decreasing, there is a rank $t_{0}$ such that for all $t \geq t_{0}$, $\lVert I_{d} - \gamma_{t} D \rVert_{\op} \leq 1- \lambda_{\min} (D) \gamma_{t}$. Then, for all $T \geq t_{0}$,
\begin{align*}
\left\| \beta_{T,0} \right\|_{\op} & \leq \prod_{t=1}^{t_{0}-1} \left( 1+ \gamma_{t}\lambda_{\max}(D) \right) \prod_{t=t_{0}}^{T} \left( 1- \gamma_{t}\lambda_{\min}(D) \right)  \\
& \leq 
\exp \left(  2\lambda_{\max} (D) \frac{c_{\gamma}C_{\rho}^{\beta}}{1 + \beta \rho - \gamma} t_{0}^{1+ \beta \rho - \gamma}   \right) \exp \left( - \lambda_{\min}(D)\frac{c_{\gamma}C_{\rho}^{\beta}}{1+\beta \rho - \gamma} T^{1 + \beta \rho - \gamma} \right).
\end{align*}
With $N_{T}$ denoting $\sum_{t=1}^{T} n_{t}$, one has $T =  \frac{N_{T}}{n}$ in  the case of the constant mini-batch size, and $T \sim \left( \frac{1+\rho}{C_{\rho}}N_{T} \right)^{\frac{1}{1+\rho}} $ for the increasing mini-batch size. Then, one has
\begin{equation}
 \label{rateinit}   \left\| \beta_{T,0} \right\|_{\op} = \left\lbrace \begin{array}{lr}
      \mathcal{O} \left( \exp \left(  - \lambda_{\min}(D) \frac{c_{\gamma}n^{\gamma - 1}}{1- \gamma}N_{T}^{1 - \gamma}  \right) \right)   &\text{if $n_{t} = n $}, \\
        \mathcal{O} \left( \exp \left(  - \lambda_{\min}(D) \frac{c_{\gamma}C_{\rho}^{ \frac{\beta   -1 + \gamma}{1+\rho}}}{1+\beta\rho-\gamma}N_{T}^{\frac{1+ \beta \rho - \gamma}{1+\rho}}  \right) \right)   &\text{if $n_{t} = \left\lfloor C_{\rho} t^{\rho} \right\rfloor $ .}
    \end{array}  \right.
\end{equation}
Then, in both cases, this term converges exponentially fast to zero.

\paragraph*{A first rate of convergence of $M_{T}$.} First, remark that 
\begin{equation}\label{maj::xi}
\mathbb{E}\left[ \left\| \xi_{t+1} \right\|^{2} |\mathcal{F}_{t} \right] \leq \mathbb{E}\left[ \left\| f_{t+1}' \right\|^{2} |\mathcal{F}_{t} \right]  \leq \frac{1}{n_{t+1}}C \left( 1+ F \left( \theta_{t} \right) - F(\theta^{*}) \right) + \left\| \nabla F \left( \theta_{t} \right) \right\|^{2} .
\end{equation}
 Then, applying \citet[Theorem 6.1]{cenac2020efficient}, one has, since $A_{t}$ converges almost surely to $A$, that
\begin{equation}
\label{firstratemartingale} \left\| M_{T} \right\|^{2} = \mathcal{O} \left( \ln (T) T^{\beta \rho - \gamma}  \right) \; \as
\end{equation}
Observe that for the constant mini-batch size, we already have the good rate of convergence for this term, but not for the increasing case. We will come back later to this term below when we find the first rate of convergence of $\theta_{T}$.

\paragraph*{A first rate of convergence of $M_{2,T}$.} As $\left\| \delta_{t} \right\| = o \left( \left\| \theta_{t} - \theta^{*} \right\| \right)$ a.s and  $\left\| A_{t} - A \right\|_{\op}$ converge almost surely to $0$, there exists a sequence of random positive variables $r_{t}$ which converges to $0$ almost surely, such that for all $t \geq t_{0}$,
\begin{align*}
  \notag \left\| R_{2,t+1} \right\| & \leq \left( 1- \gamma_{t+1} \right) \left\| R_{2,t} \right\| + \gamma_{t+1} r_{t+1}\left\| \theta_{t} - \theta^{*} \right\|   \\
  & \leq \left( 1- \gamma_{t+1} \right) \left\| R_{2,t} \right\| + 2\gamma_{t+1} r_{t+1} \left( \left\| R_{2,t} \right\|^{2} + \left\| M_{t} + R_{1,t} \right\|  \right).
\end{align*}
Then, with the help of \eqref{rateinit} and \eqref{firstratemartingale}, there exists a positive random variable $C_{1}$, such that
\begin{equation}
    \label{decR2T} \left\| R_{2,t+1} \right\|^{2}  \leq  \left( 1- \gamma_{t+1} \right) \left\| R_{2,t} \right\| + 2\gamma_{t+1} r_{t+1} \left( \left\| R_{2,t} \right\|  + C_{1} \ln( t+1)(t+1)^{\frac{\beta \rho - \gamma}{2}} \right),
\end{equation}
so that
\[
\left\| R_{2,T} \right\|^{2} = \mathcal{O} \left(  \ln (T) T^{\beta \rho - \gamma} \right) \; \as
\]
This concludes the proof for the constant mini-batch size case. For the non constant case, we need to return to the martingale term.

\paragraph*{A good rate of convergence for $M_{T}$ and $R_{2,T}$.}  Let $k_{0} = \inf \left\lbrace k , k \left(\gamma - \beta \rho \right) > \rho \right\rbrace $. Then, let us prove by induction that for any non negative integer $k \leq  k_{0}$,
\[
\left\| \theta_{T} - \theta^{*} \right\|^{2} = O \left( \ln (T)^{k}T^{-k \left( \gamma - \beta \rho \right)} \right) \; \as
\]
If $k_{0} = 0$, this is satisfied. Let us suppose from now on that $k_{0} \geq 1$ and prove this result by induction: Suppose it is true for $k-1$. Then, thanks to inequality \cref{maj::xi}, one has
\[
\mathbb{E}\left[ \left\| \xi_{t+1} \right\|^{2} |\mathcal{F}_{t} \right]  = \mathcal{O} \left( \ln (T)^{k-1}  T^{  -(k-1) \left( \gamma - \beta \rho \right)  }  \right) \; \as,
\] 
and with the help of \citet[Theorem 6.1]{cenac2020efficient}, we have
\[
\left\| M_{T} \right\|^{2} = \mathcal{O} \left( \ln (T)^{k} T^{- k \left( \gamma - \beta \rho \right) }  \right) \; \as,
\]
and $\left\| R_{2,T} \right\|^{2} = \mathcal{O} \left( \ln (T)^{k} T^{- k \left( \gamma - \beta \rho \right) }  \right) $ a.s., which concludes the induction proof.

As a particular case, one has
\[
\left\| \theta_{T} - \theta^{*} \right\|^{2} = O \left( \ln (T)^{k_{0}}T^{-  k_{0}   \left( \gamma - \beta \rho \right)} \right) \; \as,
\]
so that by definition of $k_{0}$, $\mathbb{E}[ \lVert \xi_{t+1} \rVert^{2} |\mathcal{F}_{t} ] = \mathcal{O} \left( t^{-\rho} \right)$ a.s., and we obtain with the help of \citet[Theorem 6.1]{cenac2020efficient}, that
\[
\left\| M_{T} \right\|^{2} = \mathcal{O} \left( \ln (T)  T^{-\rho - \gamma + \beta \rho }  \right) \; \as \qquad \text{and} \qquad  \left\| R_{2,T} \right\|^{2} = \mathcal{O} \left( \ln (T)  T^{-\rho - \gamma + \beta \rho }  \right) \; \as
\]
Then, since $T \sim \left( \frac{1+\rho}{C_{\rho}}N_{T} \right)^{\frac{1}{1+\rho}} $, one has
\[
\left\| \theta_{T} - \theta^{*} \right\|^{2} = O \left( \ln \left( N_{T} \right)  N_{T}^{\frac{-\rho - \gamma + \beta \rho }{1+\rho}}  \right) \; \as
\]

\subsection{Proof of \cref{theotlc,theotlc::increasing}}
Observe that one has for all $t \geq 0$,
\begin{align*}
    \theta_{t+1} - \theta^{*} = \theta_{t} - \theta^{*} - \gamma_{t+1}A H \left( \theta_{t} - \theta^{*} \right) - \gamma_{t+1} \left( A_{t} - A \right)H \left( \theta_{t} - \theta^{*} \right)   + \gamma_{t+1}A_{t} \xi_{t+1} - \gamma_{t+1}A_{t} \delta_{t},
\end{align*}
which can be written as
\begin{align}\label{decdirecte}
     \theta_{t} - \theta^{*}   = H^{-1}A^{-1}\frac{u_{t}  - u_{t+1}}{\gamma_{t+1}}   + H^{-1} A^{-1}A_{t} \xi_{t+1} - H^{-1} A^{-1}A_{t} \delta_{t} -  H^{-1}A^{-1}\left( A_{t} - A \right)H \left( \theta_{t} - \theta^{*} \right),
\end{align}
where $u_{t} = \theta_{t} - \theta^{*}$. Summing these equalities and dividing by $s_{T} = \sum_{t=0}^{T-1}n_{t+1}\ln (t+1)^{w}$, we have
\begin{align}
 \notag  \theta_{T,w} - \theta^{*} & = H^{-1}A^{-1} \frac{1}{s_{T}} \sum_{t=0}^{T-1}n_{t+1} \ln (t+1)^{w}\frac{u_{t}  - u_{t+1}}{\gamma_{t+1}}  
    + H^{-1}A^{-1} \frac{1}{s_{T}} \sum_{t=0}^{T-1}n_{t+1} \ln (t+1)^{w} A_{t} \xi_{t+1} \\
   \label{decmoy}     & - H^{-1}A^{-1} \frac{1}{s_{T}} \sum_{t=0}^{T-1}n_{t+1} \ln (t+1)^{w} A_{t}\delta_{t} -  \frac{1}{s_{T}}\sum_{t=0}^{T-1} n_{t+1}\ln(t+1)^{w}H^{-1}A^{-1}\left( A_{t} - A \right)H \left( \theta_{t} - \theta^{*} \right).
\end{align}
The rest of this proof consists in giving the rate of convergence of each term on the right-hand side of previous decomposition.

\paragraph*{Rate of convergence of $H^{-1}A^{-1} \frac{1}{s_{T}} \sum_{t=0}^{T-1}n_{t+1} \ln (t+1)^{w}A_{t} \xi_{t+1}$. } Remark that $M_{T}'= \sum_{t=0}^{T-1}n_{t+1} \ln (t+1)^{w}A_{t} \xi_{t+1}$ is a martingale term and that 
\begin{align*}
  \left\langle M_{T}' \right\rangle   & = \sum_{t=0}^{T-1} n_{t+1}^{2} \ln (t+1)^{2w}A_{t}  \mathbb{E} \left[ \xi_{t+1}\xi_{t+1}^{T} |\mathcal{F}_{t} \right]A_{t}  \\
&  = \sum_{t=0}^{T-1} n_{t+1} \ln (t+1)^{2w} A_{t} \mathbb{E}\left[  f_{t+1,i}' f_{t+1,i}'^{\top} |\mathcal{F}_{t} \right]  A_{t} \\
& - \sum_{t=0}^{T-1} n_{t+1}^{2} \ln (t+1)^{2w}  A_{t} \nabla_{\theta} F \left( \theta_{t} \right) \nabla_{\theta} F \left( \theta_{t} \right)^{T}  A_{t}.
\end{align*}
Since
\[
n_{t+1} \left\| \nabla F \left( \theta_{t} \right) \right\|^{2} = \mathcal{O} \left(\frac{\ln t}{t^{\gamma - \beta \rho }} \right) \; \as,
\]
then this converges to $0$. Next, since $\theta_{t}$ and $A_{t}$  converge to $\theta^{*}$ and $A$, we have
\[
\frac{1}{\sum_{t=0}^{T-1}n_{t+1}\ln(t+1)^{2\omega}} \left\langle M_{T}' \right\rangle\xrightarrow[T \to + \infty]{\as} A \Sigma  A .
\]
Then, with the help of a law of large numbers for martingales, we obtain that
\[
\frac{1}{s_{T}^{2}}\left\| M_{T}' \right\|^{2}  = \mathcal{O} \left( \frac{\sum_{t=0}^{T-1}n_{t+1}\ln(t+1)^{2\omega} \ln \left( \sum_{t=0}^{T-1}n_{t+1}\ln(t+1)^{2\omega} \right)}{s_{T}^{2}} \right) \; \as,
\]
which can be written as
\[
\frac{1}{s_{T}^{2}}\left\| M_{T}' \right\|^{2}  = O \left(  \frac{\ln(T+1)}{T^{\rho +1}} \right) \; \as.
\]
This, can also be written as
\[
\frac{1}{s_{T}^{2}}\left\| M_{T}' \right\|^{2} = \mathcal{O} \left( \frac{\ln(N_{T})}{N_{T}} \right) \; \as
\]
In addition, Central Limit Theorem for martingales yields,
\[
\frac{1}{\sqrt{\sum_{t=0}^{T-1}n_{t+1}\ln(t+1)^{2\omega}}}M_{T}' \xrightarrow[n\to + \infty]{\mathcal{L}}\mathcal{N}\left( 0 , A \Sigma  A \right) .
\]
Thus, as
\[
\frac{\sqrt{N_{T}}\sqrt{\sum_{t=0}^{T-1}n_{t+1}\ln(t+1)^{2\omega}}}{s_{T}} \xrightarrow[T\to + \infty]{\as} 1 ,
\]
we have
\[
\sqrt{N_{T}} \frac{1}{s_{T}}H^{-1}A^{-1} M_{T}' \xrightarrow[T\to + \infty]{\mathcal{L}} \mathcal{N}\left( 0 ,H^{-1}\Sigma H^{-1} \right) .
\]

\paragraph*{Rate of convergence of $H^{-1}A^{-1} \frac{1}{s_{T}} \sum_{t=0}^{T-1}n_{t+1} \ln (t+1)^{w} \frac{u_{t}  - u_{t+1}}{\gamma_{t+1}}$.} With the help of Abel's transformation, one have
\begin{align*}
&\frac{1}{s_{T}} \sum_{t=0}^{T-1}n_{t+1} \ln (t+1)^{w}\frac{u_{t}  - u_{t+1}}{\gamma_{t+1}} \\ & = - \frac{u_{T}n_{T}\ln(T)^{w}}{\gamma_{T}s_{T}} + \frac{u_{0}n_{1}\mathbb{1}_{\{w = 0\}}}{\gamma_{1}s_{T}} + \frac{1}{s_{T}}\sum_{t=1}^{T-1} u_{t} \left( \frac{n_{t+1}\ln(t+1)^{w}}{\gamma_{t+1}} - \frac{n_{t}\ln(t)^{w}}{\gamma_{t}}  \right).
\end{align*}
One has thanks to \cref{theo::rate,theo::rate::increasing}, we have
\[
\left\|  \frac{u_{T}n_{T}\ln(T)^{w}}{\gamma_{T}s_{T}} \right\| = \mathcal{O} \left( \frac{\sqrt{\ln T}}{T^{\frac{2+\rho - \gamma + \beta \rho}{2}}} \right) \; \as,
\]
which can be written as
\[
\left\|  \frac{u_{T}n_{T}\ln(T)^{w}}{\gamma_{T}s_{T}} \right\| = \mathcal{O} \left( \frac{\sqrt{\ln N_T}}{N_{T}^{\frac{2+\rho - \gamma + \beta \rho}{2(1+\rho)}}} \right) \; \as,
\]
which is negligible as soon as $\gamma - \beta \rho < 1$.
In addition, it is obvious that $\frac{u_{0}n_{1}\mathbb{1}_{\{w = 0\}}}{\gamma_{1}s_{T}} $ is negligible too. Furthermore, observe that 
  \[
\left| \frac{n_{t+1}\ln(t+1)^{w}}{\gamma_{t+1}} - \frac{n_{t}\ln(t)^{w}}{\gamma_{t}} \right| \leq C_{\rho} \max \left\lbrace \rho (1-\beta) + \gamma , w \right\rbrace  \max \left\lbrace t^{\rho(1-\beta) + \gamma -1} , (t+1)^{\rho(1-\beta) + \gamma -1} \right\rbrace  \ln(t+1)^{w},
  \]
which with the help of \cref{theo::rate,theo::rate::increasing} yields,
\[
\left\| \sum_{t=0}^{T-1} \left( \frac{n_{t+1}\ln(t+1)^{w}}{\gamma_{t+1}} - \frac{n_{t}\ln(t)^{w}}{\gamma_{t}} \right)  \left(   \theta_{t} - \theta^{*} \right) \right\| = \mathcal{O} \left( \ln(T)^{w+1/2} T^{\frac{\rho(1-\beta) + \gamma}{2}} \right) \; \as
\]
From this, we have
\[
\frac{1}{s_{T}}\left\| \sum_{t=0}^{T-1} \left( \frac{n_{t+1}\ln(t+1)^{w}}{\gamma_{t+1}} - \frac{n_{t}\ln(t)^{w}}{\gamma_{t}} \right)   \left( \theta_{t} - \theta^{*} \right) \right\| = \mathcal{O} \left( \sqrt{\ln(T)} T^{  \frac{ -2+ \gamma-\rho(1+\beta)}{2}} \right) \; \as,
\]
which can be written as
\begin{equation}\label{noredite}
\frac{1}{s_{T}}\left\| \sum_{t=0}^{T-1} \left( \frac{n_{t+1}\ln(t+1)^{w}}{\gamma_{t+1}} - \frac{n_{t}\ln(t)^{w}}{\gamma_{t}} \right)   \left( \theta_{t} - \theta^{*} \right) \right\| = \mathcal{O} \left( \sqrt{\ln\left( N_{T} \right)} N_{T}^{  \frac{ -  2+ \gamma-\rho(1+\beta)}{2(1+\rho)}} \right) \; \as,
\end{equation}
which is negligible as soon as $\gamma - \beta \rho < 1$.

\paragraph*{Rate of convergence of $H^{-1}A^{-1}\frac{1}{s_{T}} \sum_{t=0}^{T-1}n_{t+1} \ln (t+1)^{w}\left( A_{t} - A \right) H \left( \theta_{t} - \theta^{*} \right)$. }
Since $\left\| A_{t} - A \right\|_{\op} = \mathcal{O} \left( t^{-\nu} \right)$ a.s and with the help of \cref{theo::rate}, we have
\[
\left\| \frac{1}{s_{T}} \sum_{t=0}^{T-1}n_{t+1} \ln (t+1)^{w}\left( A_{t} - A \right) H \left( \theta_{t} - \theta^{*} \right) \right\|= \left\lbrace 
\begin{array}{lr}
\mathcal{O} \left( \frac{(\ln T)^{ \frac{1}{2}+ \mathbb{1}_{\nu + \frac{\rho(1-\beta) + \gamma}{2}=1}}}{T^{{  \nu + \frac{\rho(1-\beta) + \gamma}{2}}}} \right)  \;\as     & \text{ if } \nu + \frac{\rho(1-\beta) + \gamma}{2} \leq 1 \\
\mathcal{O} \left( \frac{1}{T^{1+\rho}} \right) \;\as    &  \text{else}
\end{array} \right.
\]
which can be written as
\[
\left\| \frac{1}{s_{T}} \sum_{t=0}^{T-1}n_{t+1} \ln (t+1)^{w}\left( A_{t} - A \right) H \left( \theta_{t} - \theta^{*} \right) \right\|= \left\lbrace  \begin{array}{lr}
\mathcal{O} \left( \frac{\ln \left( N_{T} \right)^{ \left( \frac{1}{2}+ \mathbb{1}_{\nu + \frac{\rho(1-\beta) + \gamma}{2}=1} \right)}}{N_{T}^{\frac{   2\nu + \rho(1-\beta) + \gamma}{2(1+\rho)}}} \right) \;\as     & \text{ if } \nu + \frac{\rho(1-\beta) + \gamma}{2} \leq 1 \\
\mathcal{O} \left( \frac{1}{N_{T}} \right)  \;\as   &  \text{else}
\end{array} \right.
\]

\paragraph*{Rate of convergence of $H^{-1}A^{-1} \frac{1}{s_{T}} \sum_{t=0}^{T-1}n_{t+1} \ln (t+1)^{w} A_{t}\delta_{t}$. } As $\left\| \delta_{t} \right\| \leq L_{\delta} \left\| \theta_{t} - \theta^{*} \right\|^{2}$ and with the help of \cref{theo::rate}, we have
\[
\left\| H^{-1}A^{-1} \frac{1}{s_{T}} \sum_{t=0}^{T-1}n_{t+1} \ln (t+1)^{w} A_{t}\delta_{t} \right\| = O \left(   \frac{\ln T}{T^{\rho (1-\beta  ) + \gamma}} \right)  \; \as,
\]
which can be written as
\[
\left\| H^{-1}A^{-1} \frac{1}{s_{T}} \sum_{t=0}^{T-1}n_{t+1} \ln (t+1)^{w} A_{t}\delta_{t} \right\| = O \left(   \frac{\ln N_{T}}{N_{T}^{\frac{\rho (1-\beta  ) + \gamma}{1+\rho}}} \right)  \;\as,
\]
which is negligible as soon as $\gamma > \frac{\rho (2\beta -1) +1}{2}$.

\subsection{Proof of \cref{theowithoutrate,theowithoutrate::increasing}}
Fisrt, remark that one can rewrite decomposition \cref{decdirecte} to
\[
\theta_{t} - \theta^{*} = H^{-1}A_{t}^{-1} \frac{u_{t} - u_{t+1}}{\gamma_{t+1}} + H^{-1}\xi_{t+1}  - H^{-1}\delta_{t} ,
\]
meaning that
\begin{align*}
    \theta_{T,w} - \theta^{*} &  = \frac{1}{s_{T}}   \sum_{t=0}^{T-1} n_{t+1} \ln (t+1)^{w}H^{-1}A_{t}^{-1} \frac{u_{t} - u_{t+1}}{\gamma_{t+1}} +\frac{1}{S_{T}}   \sum_{t=0}^{T-1} n_{t+1} \ln (t+1)^{w} H^{-1}\xi_{t+1} \\
    &  - \frac{1}{s_{T}}   \sum_{t=0}^{T-1} n_{t+1} \ln (t+1)^{w}H^{-1}\delta_{t} .
\end{align*}
Analogously to the proof of \cref{theotlc}, one can easily check that
\[
 \left\| \frac{1}{S_{T}}   \sum_{t=0}^{T-1} n_{t+1} \ln (t+1)^{w} H^{-1}\xi_{t+1} \right\|^{2} = \mathcal{O} \left( \frac{\ln N_{T}}{N_{T}} \right)  \;\as,
\]
and
\[
 \quad \quad \sqrt{N_{T}}\frac{1}{s_{T}}   \sum_{t=0}^{T-1} n_{t+1} \ln (t+1)^{w} H^{-1}\xi_{t+1} \xrightarrow[T\to + \infty]{\mathcal{L}} \mathcal{N} \left( 0 , H^{-1}\Sigma H^{-1} \right) .
\]
In the same way, we have
\[
\left\| \frac{1}{s_{T}}   \sum_{t=0}^{T-1} n_{t+1} \ln (t+1)^{w}H^{-1}\delta_{t}  \right\| = O \left(   \frac{\ln N_{T}}{N_{T}^{\frac{\rho (1-\beta  ) + \gamma}{1+\rho}}} \right)  \;\as
\]
In addition, note that 
\begin{align*}
     \frac{1}{s_{T}} \sum_{t=0}^{T-1} n_{t+1}    \ln (t+1)^{w}H^{-1}A_{t}^{-1} \frac{u_{t} - u_{t+1}}{\gamma_{t+1}} = &  \frac{1}{s_{T}}   \sum_{t=0}^{T-1} n_{t+1} \ln (t+1)^{w}H^{-1} \frac{A_{t}^{-1}u_{t} - A_{t+1}^{-1}u_{t+1}}{\gamma_{t+1}}  \\ & +  \frac{1}{s_{T}}   \sum_{t=0}^{T-1} n_{t+1} \ln (t+1)^{w}H^{-1}\left( A_{t +1}^{-1} - A_{t}^{-1} \right) \frac{ u_{t+1}}{\gamma_{t+1}} .
\end{align*}
With the help of Abel's transformation and since $A_{t}$ converges almost surely to the positive matrix $A$ (\cref{ass::4}), following the lines of the proof for \cref{theotlc} (e.g., see \cref{noredite}), one can show that
\[
\left\| \frac{1}{s_{T}}   \sum_{t=0}^{T-1} n_{t+1} \ln (t+1)^{w}H^{-1} \frac{A_{t}^{-1}u_{t} - A_{t+1}^{-1}u_{t+1}}{\gamma_{t+1}} \right\| = O \left( \frac{\sqrt{ \ln\left( N_{T} \right)}}{ N_{T}^{  \frac{ 2+ \gamma-\rho(1+\beta)}{2(1+\rho)}}} \right) \; \as
\]
In addition, since $\left\| \theta_{t} - \theta^{*} \right\|^{2} = \mathcal{O} \left( \frac{\ln(t)}{t^{\gamma - \beta \rho + \rho}} \right)$ a.s., with $E_{t}$ denoting the event $\lbrace \lVert \theta_{t} - \theta \rVert^{2} \leq \frac{(\ln(t))^{1+\delta}}{t^{\gamma +  \rho (1- \beta)}} , \lVert \theta_{t,w} - \theta \rVert^{2} \leq \frac{(\ln(t))^{1+\delta}}{t^{\gamma + \rho(1-  \beta )}}\rbrace $,  $\mathbb{1}_{\{E_{t}^{C}\}}$ converges almost surely to $0$, then, we have
\[
 \frac{1}{s_{T}}   \sum_{t=0}^{T-1} n_{t+1} \ln (t+1)^{w}H^{-1}\left\| A_{t+1}^{-1} - A_{t}^{-1} \right\|_{\op} \frac{ \left\| u_{t+1} \right\|}{\gamma_{t+1}}  \mathbb{1}_{\{E_{t+1}^{C}\}} = \mathcal{O} \left( \frac{1}{N_{T}} \right) \; \as,
\]
and
\begin{align*}
\frac{1}{s_{T}}   \sum_{t=0}^{T-1} n_{t+1} \ln (t+1)^{w}H^{-1} & \left\| A_{t+1}^{-1} - A_{t}^{-1} \right\|_{\op} \frac{ \left\| u_{t+1} \right\|}{\gamma_{t+1}} \mathbb{1}_{\{E_{t+1}\}}  \\
& \leq \frac{1}{s_{T}}   \sum_{t=0}^{T-1} n_{t+1} \ln (t+1)^{w + 1/2 + \delta}H^{-1}\left\| A_{t+1}^{-1} - A_{t}^{-1} \right\|_{\op} c_{\gamma}^{-1} (t+1)^{\frac{\gamma -  \rho (\beta +1)}{2}}.
\end{align*}
At last, one can conclude the proof with the help of equality \cref{equalitystrange::increasing}.

\subsection{Proof of \cref{theo::newton::direct}}
Observe that the convergence of $\theta_{T}$ is obtained with the same calculus as in the proof of \cref{theo::ps}.
Remark that decomposition \cref{decdelta} can now be written as
\[
\theta_{t+1} - \theta^{*} =  \left( 1 - \frac{n_{t+1}}{N_{t+1}} \right)  \left( \theta_{t} - \theta^{*} \right)  - \frac{n_{t+1}}{N_{t+1}}\left( \overline{H}_{t}^{-1} - H^{-1} \right) H \left( \theta_{t} - \theta^{*} \right) + \frac{n_{t+1}}{N_{t+1}}  \overline{H}_{t}^{-1} \xi_{t+1} - \frac{n_{t+1}}{N_{t+1}}  \overline{H}_{t}^{-1} \delta_{t} .
\]
Then, with the help of induction, one has
\begin{align*}
\theta_{T} - \theta^{*} = &  \frac{1}{N_{T}}\left( \theta_{0} - \theta^{*} \right) \underbrace{ - \frac{1}{N_{T}}\sum_{t=0}^{T-1} n_{t+1} \left(  \overline{H}_{t}^{-1} - H^{-1} \right)H \left( \theta_{t} - \theta^{*} \right)  - \frac{1}{N_{T}}\sum_{t=0}^{T-1}n_{t+1} \overline{H}_{t}^{-1} \delta_{t}}_{=: \Delta_{T}} \\ & + \frac{1}{N_{T}}\underbrace{\sum_{t=0}^{T-1} n_{t+1} \overline{H}_{t}^{-1}\xi_{t+1}}_{=: M_{T}} .
\end{align*}

\paragraph{Convergence of the martingale term $M_{T}$.} Observe that $M_{T}$ is a martingale term and that
\begin{align*}
\left\langle M \right\rangle_{T} = \sum_{t=0}^{T-1} n_{t+1}^{2}  \overline{H}_{t}^{-1}\mathbb{E} \left[ \xi_{t+1}\xi_{t+1}^{T} \right] A_{t} = & \sum_{t=0}^{T-1} n_{t+1}  \overline{H}_{t}^{-1} \mathbb{E}\left[ \nabla_{\theta} f \left( X_{t+1} , \theta_{t} \right)\nabla_{\theta} f \left( X_{t+1}, \theta_{t} \right)^{T} |\mathcal{F}_{t} \right]  \overline{H}_{t}^{-1} \\
&  - \sum_{t=0}^{T-1}n_{t+1} \overline{H}_{t}^{-1}  \nabla F \left( \theta_{t} \right) \nabla F \left( \theta_{t} \right)^{T} \overline{H}_{t}^{-1}.
\end{align*}
Then, since $\theta_{t}$ and $ \overline{H}_{t}^{-1}$ converge  almost surely to $\theta^{*}$ and $H^{-1}$ and by continuity (\cref{ass::1}), one obtain that
\[
\frac{1}{N_{T}}\left\langle M \right\rangle_{T}  \xrightarrow[T\to + \infty]{\as}  H^{-1}\Sigma H^{-1} . 
\]
Thus, with the help of a law of large numbers for martingales, we have
\[
 \left\| \frac{1}{N_{T}} M_{T} \right\|^{2} = \mathcal{O} \left( \frac{\ln N_{T}}{N_{T}} \right) \; \as,
\]
and with the help of Central Limit Theorem for martingales,
\[
\frac{1}{\sqrt{N_{T}}}M_{T}\xrightarrow[n\to + \infty]{\mathcal{L}} \mathcal{N}\left( 0 ,   H^{-1}\Sigma H^{-1} \right)   .
\]

\paragraph{Convergence of the rest terms.} Since $\overline{H}_{t}^{-1}$ converges to $H^{-1}$ and $\left\| \delta_{t} \right\| = o \left( \left\| \theta_{t} - \theta^{*} \right\| \right)$ a.s., there is a sequence of positive random variables $(r_{t}')$ converging to $0$, such that
\begin{align*}
\left\| \Delta_{T+1} \right\| &   \leq \left( 1- \frac{n_{T+1}}{N_{T+1}} \right) \left\| \Delta_{T} \right\| + \frac{n_{T+1}}{N_{T+1}} r_{T}' \left\| \theta_{T} - \theta^{*} \right\| \\
& \leq \left( 1- \frac{n_{T+1}}{N_{T+1}} \right) \left\| \Delta_{T} \right\| + \frac{n_{T+1}}{N_{T+1}} r_{T}' \left\| \frac{1}{N_{T}} \left( \theta_{0} - \theta^{*} \right) + \frac{1}{N_{T}} M_{T} + \Delta_{T} \right\|  .
\end{align*}
Then, there is a positive random variable $C_{M}$, such that 
\[
\left\| \Delta_{T+1} \right\| \leq \left( 1- \frac{n_{T+1}}{N_{T+1}} \right) \left\| \Delta_{T} \right\| + \frac{n_{T+1}}{n_{T+1}}r_{T}' \left( \left\| \Delta_{T} \right\| + C_{M} \frac{\sqrt{\ln T}}{T^{\frac{1+\rho}{2}}} \right),
\]
which can also be written for any $c\in (0,1)$ as
\[
\left\| \Delta_{T+1} \right\| \leq \left( 1- c\frac{n_{T+1}}{N_{T+1}} \right) \left\| \Delta_{T} \right\| +   \frac{n_{T+1}}{n_{T+1}}   C_{M} \frac{\sqrt{\ln (T+1)}}{(T+1)^{\frac{1+\rho}{2}}}  + r_{T}'' ,
\]
with $r_{T}'' =   \frac{n_{T+1}}{n_{T+1}} r_{T}' \left( \left\| \Delta_{T} \right\| + C_{M} \frac{\sqrt{\ln (T+1)}}{(T+1)^{\frac{1+\rho}{2}}} \right) \mathbb{1}_{r_{T}' > c}$. Then, with the help of an induction, one has
\begin{align*}
\left\| \Delta_{T} \right\| \leq \tilde{\beta}_{T,0} \left\| \Delta_{0} \right\|  + \sum_{t=0}^{T-1} \tilde{\beta}_{T,t+1} \frac{n_{t+1}}{N_{t+1}} C_{M} \frac{\sqrt{\ln (t+1)}}{(t+1)^{\frac{1+\rho}{2}}} + \sum_{t=0}^{T-1} \tilde{\beta}_{T,t+1} r_{t}'',
\end{align*}
with $\tilde{\beta}_{T,t} = \prod_{j=t+1}^{T}  \left( 1- c\frac{n_{j}}{N_{j}} \right)$ and $\beta_{T,T} = 1$.
In addition, since for any $t$, one has $N_{t} \leq \frac{C_{\rho}}{1+\rho} \left( (t+1)^{1+\rho} -1 \right)$, one has for any $t \leq T$,
\begin{align*}
\tilde{\beta}_{T,t} \leq & \exp \left( - c \sum_{j=t+1}^{T} \frac{n_{j}}{N_{j}} \right)  \leq \exp \left( - c(1+\rho) \sum_{j=t+1}^{T}\frac{j^{\rho}}{(j+1)^{1+\rho} } \right) \\ \leq & \exp \left( -  c(1+\rho) \left( \frac{t+1}{t+2} \right)^{\rho} \sum_{j=t+1}^{T}\frac{1}{j+1} \right)  \leq \left( \frac{ t+1 }{T+1} \right)^{c_{t}}  ,
\end{align*}
with $c_{t} = c(1+\rho) \left( \frac{t+1}{t+2} \right)^{\rho} \geq c(1+\rho)2^{-\rho}$. Taking $1>c > 2^{\rho -1} $ and denoting $c_{\rho} = c2^{-\rho} > 1/2$, one has
\[
\tilde{\beta}_{T,t}  \leq \left( \frac{ t+1 }{T+1} \right)^{c_{\rho}(1+\rho)}  .
\]
Then, as a particular case, $\tilde{\beta}_{T,0} \leq \frac{1}{(T+1)^{c_{\rho}(1+\rho)}}$ and this term is so negligible. In addition, since 
\[
\sum_{t=0}^{T-1} \tilde{\beta}_{T,t+1} r_{t}'' =\tilde{\beta}_{T,0} \sum_{t=0}^{T-1} \tilde{\beta}_{t+1,0}^{-1} r_{t}'',
\]
and since $ \mathbb{1}_{\{r_{T}' > c\}}$ converges almost surely to $0$, one has
\begin{align*}
\sum_{t=0}^{T-1} \tilde{\beta}_{T,t+1} r_{t}'' = \mathcal{O} \left( \tilde{\beta}_{T,0}  \right)  = \mathcal{O} \left(    \frac{1}{(T+1)^{c_{\rho}(1+\rho)}} \right) \; \as,
\end{align*}
and this term is so negligible as $c_{\rho}> 1/2$.
Finally, 
\begin{align*}
\sum_{t=0}^{T-1} \tilde{\beta}_{T,t+1} \frac{n_{t+1}}{N_{t+1}} C_{M} \frac{\sqrt{\ln (t+1)}}{(t+1)^{\frac{1+\rho}{2}}} \leq \sum_{t=0}^{T-1} \left( \frac{ t+1 }{T+1} \right)^{c_{\rho}(1+\rho)} \frac{n_{t+1}}{N_{t+1}}C_{M}\frac{\sqrt{\ln (t+1)}}{(t+1)^{\frac{1+\rho}{2}}} = \mathcal{O} \left(  \frac{\sqrt{\ln T}}{T^{\frac{1+\rho}{2}}} \right) \;\as,
\end{align*}
leading to $\left\| \Delta_{T} \right\| = \mathcal{O} \left( \sqrt{\frac{\ln T}{T^{1+\rho}}} \right)$ a.s., and 
\[
\left\| \theta_{T} - \theta^{*} \right\|^{2} = O \left( \frac{\ln T}{T^{1+\rho}} \right) \;\as,  \qquad \text{and} \qquad \left\| \theta_{T} - \theta^{*} \right\|^{2} = O \left( \frac{\ln N_T}{N_T} \right) \; \as
\]

\paragraph{Asymptotic efficiency.}
In order to get the asymptotic normality, we now have to give a better rate of convergence of $\left\| \Delta_{T} \right\|$. First, since $ \bar{H}_{t}^{-1}$ converges to $H^{-1}$, $\left\| \delta_{t} \right\| \leq  L_{\delta} \left\| \theta_{t}  - \theta^{*} \right\|^{2}$, and with the help of the rate of convergence of $\theta_{t}$, one has
\[
\frac{1}{N_T} \left\| \sum_{t=0}^{T-1} n_{t+1}  \bar{H}_{t}^{-1} \delta_{t} \right\| \leq \frac{L_{\delta}}{N_T}\sum_{t=0}^{T-1} n_{t+1} \left\| \bar{H}_{t}^{-1} \right\|_{\op} \left\| \theta_{t} - \theta^{*} \right\|^{2} = \mathcal{O} \left( \frac{(\ln T)^{2}}{T^{1+\rho}} \right) \; \as,
\]
which is a negligible term. In addition, since $\lVert  \bar{H}_{t}^{-1} - H^{-1} \rVert_{\op} = \mathcal{O} \left( t^{-\nu} \right)$ a.s., one has
\begin{align*}
\frac{1}{N_T}\left\| \sum_{t=0}^{T-1} n_{t+1} \left(  \bar{H}_{t}^{-1}  - H^{-1} \right) H \left( \theta_{t} - \theta^{*} \right) \right\| & \leq \frac{1}{N_T}\left\| H \right\|_{\op} \sum_{t=0}^{T-1} n_{t+1}\left\|  \bar{H}_{t}^{-1}  - H^{-1} \right\|_{\op} \left\| \theta_{t} - \theta^{*} \right\|  \\
& = \mathcal{O} \left(  \frac{\ln(T)^{1/2 + \mathbb{1}_{\{(1+\rho)/2 + \nu = 1\}}}}{T^{ \rho + \min \lbrace 1,(1-\rho)/2+ \nu \rbrace}} \right) \;\as
\end{align*}
Hence, as $\nu> 0$, this term is negligible, which thereby concludes the proof.

\subsection{Proof of \cref{theo::newton::nd,theo::newton::nd::increasing}}
Let us first check that the assumptions on the learning rate (step-sequence) are satisfied: First, since for all $t \geq 1$ and $i = 1 , \ldots ,n_{t}$,
\[
\frac{N_{Z,t,i}}{N_{t}} \xrightarrow[t\to + \infty]{\as} p,
\]
we can observe that\footnote{E.g., see \citet{godichon2024recursive,bercu2021stochastic} for more details.}
\[
\frac{d (1-\iota)}{p\ln (t+1)^{w'}N_{t}^{1-\iota}}\sum_{t'=1}^{t} \ln (t'+1)^{w'}\sum_{i=1}^{n_{t'}} Z_{t',i}\iota_{t',i} e_{t',i}e_{t',i}^{T} \xrightarrow[t\to + \infty]{\as}I_{d} ,
\]
such that
\[
\frac{ 1+ \sum_{t'=1}^{t} \ln (t'+1)^{w'} \sum_{i=1}^{n_{t}'} Z_{t',i}}{\ln (t+1)^{w'}N_{t}} \xrightarrow[t\to + \infty]{\as} p,
\]
Next, by definition of $\iota$, we have
\[
\lambda_{\max} \left( \bar{H}_{t,w'}^{-1} \right) = \mathcal{O} \left( t^{\iota (1+ \rho)} \right) \; \as, \qquad  \text{and} \qquad \sum_{t \geq 1} \frac{\gamma_{t}^{2}}{n_{t}}\lambda_{\max} \left( \bar{H}_{t-1,w'}^{-1} \right)^{2} < + \infty \; \as
\]
In addition, with $N_{Z,T} := 1+  \sum_{t=1}^{T} \ln (t+1)^{w'} \sum_{i=1}^{n_{t'}}Z_{t,i}$, one has
\begin{align*}
\frac{1}{N_{Z,T}} \sum_{t=1}^{T} \ln (t+1)^{w'} \sum_{i=1}^{n_{t}} Z_{t ,i} \alpha_{t,i} \Phi_{t,i}\Phi_{t,i}^{\top}  = & \frac{1}{N_{Z,T}} \sum_{t=1}^{T} \ln (t+1)^{w'}\nabla_{\theta}^{2}F \left( \theta_{t-1} \right) \sum_{i=1}^{n_{t}}Z_{t,i} \\
& + \frac{1}{N_{Z,T}} \sum_{t=1}^{T} \ln (t+1)^{w'}   \xi_{Z,t},
\end{align*}
where $\xi_{Z,t} :=  \sum_{i=1}^{n_{t}} Z_{t ,i}\alpha_{t,i} \Phi_{t,i}\Phi_{t,i}^{\top} -   \sum_{i=1}^{n_{t}} Z_{t,i} \nabla_{\theta}^{2}F \left( \theta_{t-1} \right)$ is a sequence of martingale differences for the  filtration $\mathcal{F}_{t-1}' = \sigma \left( X_{1,1} , \ldots ,X_{t-1,n_{t-1}},Z_{t,1} , \ldots , Z_{t,n_{t}} \right)$. Thus,
\[
\mathbb{E}\left[ \left\| \xi_{Z,t} \right\|_{F}^{\eta '} |\mathcal{F}_{t-1}' \right] \leq  2^{\eta ' -1}  \left( \sum_{i=1}^{n_{t}}   Z_{t ,i}  \mathbb{E} \left[  \left\| 	\alpha_{t,i} \Phi_{t,i}\Phi_{t,i}^{\top} \right\|^{\eta '} |\mathcal{F}_{t-1} \right] ^{\frac{1}{\eta '}} \right)^{\eta '} \leq 2^{\eta ' -1} C_{\eta '}^{\eta '}  \left( \sum_{i=1}^{n_t} Z_{t,i} \right)^{\eta '},
\]
and with the help of a law of large numbers for martingales, one has
\[
\left\|  \sum_{t=1}^{T} \ln (t+1)^{w'}   \xi_{Z,t} \right\|_{F} = o \left( N_{Z,T} \right) \; \as
\]
Since for all $\theta \in \mathbb{R}^{d}$,  $\left\| \nabla_{\theta}^{2}F(\theta) \right\|_{\op} \leq L_{\nabla G}$,
\[
\left\|  \frac{1}{N_{Z,T}} \sum_{t=1}^{T} \ln (t+1)^{w'}\nabla_{\theta}^{2}F \left( \theta_{t-1} \right) \sum_{i=1}^{n_{t}}Z_{t,i} \right\|_{\op} \leq L_{\nabla F} .
\]
Then,  $\lambda_{\max} \left( \bar{H}_{t,w'} \right) = \mathcal{O} (1) $ a.s., such that
\[
\sum_{t\geq 1} \gamma_{t} \lambda_{\min} \left( \bar{H}_{t-1,w'}^{-1}  \right)  = + \infty \; \as,  \qquad \text{and} \qquad \frac{\lambda_{\max} \left( \bar{H}_{t,w'} ^{-1} \right)^{2}\gamma_{t+1}}{\lambda_{\min} \left( \bar{H}_{t,w'} ^{-1} \right)} = \mathcal{O} \left( t^{2\iota (1+ \rho) - 1  } \right) \;\as,
\]
and the conditions in \cref{cond::step::increasing} are satisfied as soon as $i < \frac{1-\rho}{2(1+\rho)}$. Then, according to \cref{theo::newton::direct}, $\theta_{T}$ converges almost surely to $\theta^{*}$. By continuity, this implies that 
\[
 \frac{1}{N_{Z,T}} \sum_{t=1}^{T} \ln (t+1)^{w'}\nabla^{2}F \left( \theta_{t-1} \right) \sum_{i=1}^{n_{t}}Z_{t,i} \xrightarrow[T\to + \infty]{\as} \nabla_{\theta}^{2}F(\theta^{*})
\]
meaning that $\bar{H}_{T,w'}$ and $\bar{H}_{T,w'}^{-1}$ converge almost surely to $H$ and $H^{-1}$. Then, thanks to \cref{theo::newton::direct}, one has that
\[
\left\| \theta_{T,w} - \theta^{*} \right\|^{2} = \mathcal{O} \left( \frac{\ln N_{T}}{N_{T}} \right) \; \as,
\]
and since the Hessian is locally Lipschitz, $\left\| \nabla_{\theta}^{2} F \left( \theta_{t} \right) - H \right\|_{\op} = \mathcal{O} \left(  \sqrt{\ln N_{t}}N_{t}^{- 1/2 } \right) $ a.s., and
\[
\left\|  \frac{1}{N_{Z,T}} \sum_{t=1}^{T} \ln (t+1)^{w'}\nabla^{2}F \left( \theta_{t-1} \right) \sum_{i=1}^{n_{t}}Z_{t,i} - H \right\|_{\op}  = \mathcal{O} \left(    \frac{1}{\ln(T+1)^{w'}N_{T}} \sum_{t=1}^{T} \ln (t+1)^{w'+1/2} n_{t} t^{-\frac{1+\rho}{2}} \right) \; \as,
\]
one has
\[
\left\|  \frac{1}{N_{Z,T}} \sum_{t=1}^{T} \ln (t+1)^{w'}\nabla^{2}F \left( \theta_{t-1} \right) \sum_{i=1}^{n_{t}}Z_{t,i} - H \right\|_{\op}  = \mathcal{O} \left( \frac{\sqrt{\ln(N_{T})}}{\sqrt{N_{T}}}  \right) \; \as
\]
In addition, since $\eta ' \geq 2$ and
\[
\mathbb{E}\left[ \left\| \xi_{Z,t} \right\|_{F}^{2} |\mathcal{F}_{t-1} \right] \leq  \sum_{i=1}^{n_{t}} Z_{t,i}^{2} \mathbb{E}\left[ \left\| a \left( X_{t,i} , \theta_{t-1} \right) \Phi_{t,i}\Phi_{t,i}^{T} \right\|_{F}^{2} |\mathcal{F}_{t-1}' \right] \leq n_{t}C_{\eta'}^{\frac{2}{\eta '}} ,
\]
one has, with the help of a law of large numbers for martingales, that for all $\delta > 0$,
\[
\left\|  \frac{1}{N_{Z,T}} \sum_{t=1}^{T} \ln (t+1)^{w'}   \xi_{Z,t} \right\|_{F}^{2} = O \left( \frac{(\ln N_{T})^{1+\delta}}{N_{T}} \right) \; \as
\]
Meaning, that
\[
\left\| \frac{1}{N_{Z,T}}\sum_{t=1}^{T} \ln (t+1)^{w'} \sum_{i=1}^{n_{t}} Z_{t,i}\iota_{t,i} e_{t,i}e_{t,i}^{T} \right\|_{\op} = \mathcal{O} \left( \frac{1}{T^{\iota(1+\rho)}} \right) \; \as,
\]
and by definition of $\iota$, one has
\[
\left\| \bar{H}_{T,w'} - H \right\|^{2} = \mathcal{O} \left( \frac{1}{N_{T}^{ 2\iota }} \right) \; \as
\]
Then, with the help of \cref{theo::newton::direct}, one has
\[
\sqrt{N_{T}} \left( \theta_{T} - \theta^{*} \right) \xrightarrow[T\to + \infty]{\mathcal{L}}\mathcal{N}\left( 0 , H^{-1}\Sigma H^{-1} \right) .
\]

\subsection{Proof of \cref{theo::newton::wasn,theo::newton::wasn::inceasing}}
As in the proof of \cref{theo::newton::nd,theo::newton::nd::increasing}, one can easily check that the conditions in \cref{cond::step::increasing} are satisfied, such that \cref{theo::ps::increasing} hold, i.e., $\theta_{T}$ and $\theta_{T,w}$ converges almost surely to $\theta^{*}$. In a same way, as in the proof of \cref{theo::newton::nd}, one can easily get the consistency of $\bar{S}_{T,w'}$, leading with the help of \cref{theo::rate::increasing} to
\[
\left\| \theta_{T}  - \theta^{*} \right\|^{2} = \mathcal{O} \left( \frac{\ln (N_{T})}{N_{T}^{\frac{\gamma + \rho(1-\beta)}{1+\rho}}} \right) \;\as, \qquad \text{and} \qquad \left\| \theta_{T,w}  - \theta^{*} \right\|^{2} = \mathcal{O} \left( \frac{\ln N_{T}}{N_{T}^{\frac{\gamma + \rho(1-\beta)}{1+\rho}}} \right) \; \as
\]
In order to conclude the proof, we will now check that equality \cref{equalitystrange} is satisfied, i.e., that
\[
     \frac{1}{s_{T}}   \sum_{t=0}^{T-1} n_{t+1} \ln (t+1)^{w + 1/2 + \delta} \left( r_{t+1} + r_{t+1}' \right)  t^{\frac{\gamma -  \rho (\beta +1)}{2}} = \mathcal{O} \left( \frac{1}{T^{(1+ \rho )v'}   } \right) \; \as,
\]
with
\[
r_{t+1}'  = \frac{\ln (t+1)^{w'}}{N_{Z,t+1}}\sum_{i=1}^{n_{t+1}} Z_{t+1,i} \iota_{t+1,i} \quad \text{and} \quad r_{t+1} = \frac{\ln (t+1)^{w'}}{N_{Z,t+1}} \sum_{i=1}^{n_{t+1}}Z_{t+1,i}\left\| \alpha_{t+1,i} \Phi_{t+1,i} \right\|.
\]
First, since $\sum_{i=1}^{n_{t+1}} \iota_{t+1,i}  = \mathcal{O} \left( t^{-\iota(1+\rho) + \rho} \right)$, and since $\iota < \frac{\gamma - \rho \beta}{2(1+\rho)}  $, one has
\[
 \frac{1}{s_{T}}   \sum_{t=0}^{T-1} n_{t+1} \ln (t+1)^{w + 1/2 + \delta}   r_{t+1}'  t^{\frac{\gamma -  \rho (\beta +1)}{2}} = \mathcal{O} \left( \frac{\ln (T+1)^{  3/2 + \delta}}{T^{  \iota(1+\rho) +1 + \frac{\beta (1+\rho) - \gamma}{2}  } } \right) \; \as,
\]
and since $\gamma - \beta \rho < 1$, it comes that $\iota (1 + \rho)+1 + \frac{\beta (1+\rho) - \gamma}{2} > \frac{1+\rho}{2}$. Considering the sequence of martingale differences $\Xi_{Z,t+1} =\sum_{i=1}^{n_{t+1}}Z_{t+1,i} \alpha_{t+1,i} \left\| \Phi_{t+1,i} \right\|^{2}  - \sum_{i=1}^{n_{t+1}} Z_{t+1,i} \mathbb{E}[ \alpha_{t+1,i} \left\| \Phi_{t+1,i} \right\|^{2}  |\mathcal{F}_{t-1}' ]$, one has
\begin{align*}
  \sum_{t=0}^{T-1} n_{t+1} \ln (t+1)^{w + 1/2 + \delta}   r_{t+1}   t^{\frac{\gamma -  \rho (\beta +1)}{2}} \leq & C_{\eta'}^{\frac{1}{\eta '}} \sum_{t=0}^{T-1}\frac{ n_{t+1} \ln (t+1)^{w +w'+ 1/2 + \delta}     t^{\frac{\gamma -  \rho (\beta +1)}{2}}  }{N_{Z,t+1}}  \sum_{i=1}^{n_{t+1}} Z_{t+1,i} \\
  & + \sum_{t=0}^{T-1}\frac{ n_{t+1} \ln (t+1)^{w +w'+ 1/2 + \delta}     t^{\frac{\gamma -  \rho (\beta +1)}{2}}  }{N_{Z,t+1}}\Xi_{Z,t+1}.
\end{align*}
Furthermore,
\begin{align*}
\frac{1}{s_{T}} \sum_{t=0}^{T-1}\frac{ n_{t+1} \ln (t+1)^{w +w'+ 1/2 + \delta}     t^{\frac{\gamma -  \rho (\beta +1)}{2}}  }{N_{Z,t+1}}  \sum_{i=1}^{n_{t+1}} Z_{t+1,i} = \mathcal{O} \left( \frac{\ln(T+1)^{3/2 + \delta}}{T^{\min \left\lbrace 1 + \frac{ \rho (\beta +1) - \gamma }{2} , 1+ \rho \right\rbrace}} \right) \; \as,
\end{align*}
and since $\gamma - \beta \rho < 1$, one has that $1 + \frac{ \rho (\beta +1) - \gamma }{2}> \frac{1+\rho}{2}$. In addition, since
\[
\mathbb{E}\left[ \left\| \Xi_{Z,t+1} \right\|_{F}^{\eta '} \right] \leq \left( \sum_{i=1}^{n_{t+1}}Z_{t+1,i} \right)^{\eta '} C_{\eta '}^{\eta '},
\]
and with the help of a law of large numbers for martingales,
\begin{align*}
&\left| \sum_{t=0}^{T-1} \frac{ n_{t+1} \ln (t+1)^{w +w'+ 1/2 + \delta}     t^{\frac{\gamma -  \rho (\beta +1)}{2}}  }{N_{Z,t+1}}\Xi_{Z,t+1} \right| \\ & = o \left( \sum_{t=0}^{T-1}\frac{ n_{t+1} \ln (t+1)^{w +w'+ 1/2 + \delta}     t^{\frac{\gamma -  \rho (\beta +1)}{2}}  }{N_{Z,t+1}} \sum_{i=1}^{n_{t+1}} Z_{t+1,i} \right) \; \as,
\end{align*}
such that 
\begin{align*}
 \frac{1}{\sum_{t=0}^{T-1}n_{t+1} \ln (t+1)^{w}} & \left| \sum_{t=0}^{T-1}\frac{ n_{t+1} \ln (t+1)^{w +w'+ 1/2 + \delta}     t^{\frac{\gamma -  \rho (\beta +1)}{2}}  }{N_{Z,t+1}}\Xi_{Z,t+1} \right| \\ & = o \left( \frac{\ln(T+1)^{3/2 + \delta}}{T^{\min \left\lbrace 1 + \frac{ \rho (\beta +1) - \gamma }{2} , 1+ \rho \right\rbrace}} \right) \; \as,
\end{align*}
which concludes the proof.

\subsection{Proof of \cref{theo::adagrad,theo::adagrad::increasing}}
First, since the conditions in \cref{cond::step} (or in \cref{cond::step::increasing}) are satisfied, one has that $\theta_{t} $ and $\theta_{t,w} $ converge almost surely to $\theta^{*}$. Let us now prove that it implies the convergence of $G_{t}$.

\paragraph*{Convergence of $G_{t}$.} For all coordinate $j$, let us now consider
\[
\tilde{G}_{T}^{(j)} :=\frac{1}{N_{T}} \sum_{t=1}^{T} \sum_{i=1}^{n_{t}} \left( \frac{\partial}{\partial j} f \left( \theta_{t-1,w};\xi_{t,i}   \right) \right)^{2}.
\]
Then, denoting 
\[
V_{j} := \mathbb{E} \left[ \left(\frac{\partial}{ \partial j}  f\left( \theta^{*};\xi \right) \right)^{2}\right] =  \mathbb{V} \left[ \frac{\partial}{\partial j} f \left( \theta^{*};\xi \right) \right] ,
\]
 one has
\begin{align*}
  \tilde{G}_{T}^{(j)} - V_{j}  = \frac{1}{N_{T}}\sum_{t=1}^{T} n_{t}  \left( \mathbb{E} \left[ \left(\frac{\partial}{ \partial j}   f \left( \theta_{t-1,w} ; \xi_{t,1} \right) \right)^{2} |\mathcal{F}_{t-1}\right] - V_{j} \right) + \frac{1}{N_{T}}\sum_{t=1}^{T} \Xi_{t} 
\end{align*}
where $\Xi_{t} = \sum_{i=1}^{n_t} \left(\frac{\partial}{ \partial j}   f \left(  \theta_{t-1,w}  ;X_{t,i}\right) \right)^{2}  -  n_{t} \mathbb{E} \left[ \left(\frac{\partial}{ \partial j}   f \left(  \theta_{t-1,w} ,\xi \right) \right)^{2}\right]$ is a martingale difference. Then, thanks to \cref{momenteta::increasing} coupled with \citet[Proposition 1.III.19]{duflo2013random}, we have
\[
\frac{1}{N_{T}} \sum_{t=1}^{T}\Xi_{t} \xrightarrow[n\to + \infty]{\as} 0 .
\] 
In addition, since the functional $\theta \longmapsto \mathbb{E} \left[ \nabla_{\theta} f(\theta;\xi) \nabla_{\theta}f(\theta;\xi)^{\top} \right]$ is continuous at $\theta^{*}$, one has for all $j$
\[
\frac{1}{N_{T}}\sum_{t=1}^{T} n_{t}  \left( \mathbb{E} \left[ \left(\frac{\partial}{ \partial j}   f \left( \theta_{t-1,w},\xi_{t,1}  \right) \right)^{2} |\mathcal{F}_{t-1}\right] - V_{j} \right) \xrightarrow[T\to + \infty]{\as} 0 ,
\]
such that, for all $j$,
\[
\tilde{G}_{T}^{(j)} = \frac{1}{N_{T}} \sum_{t=1}^{T} \sum_{i=1}^{n_{t}} \left( \frac{\partial}{\partial j} f  \left(  \theta_{t-1,w};\xi_{t,i}  \right) \right)^{2} \xrightarrow[T \to + \infty]{a.s}  \mathbb{V} \left[ \frac{\partial}{\partial j} f \left(  \theta^{*};\xi \right) \right] > 0.
\]
Then, $G_{t }$ converges almost surely to the diagonal matrix $G$, whose diagonal elements are given by  $G^{(j)} = \mathbb{V} \left[ \frac{\partial}{\partial j} f(\theta^{*};\xi) \right]^{-1/2}$.

\paragraph{Rate of convergence of $\theta_{T}$.}
With the help of \cref{theo::rate::increasing}, one has that
\[
\left\| \theta_{T}  - \theta^{*} \right\|^{2} = \mathcal{O}  \left( \frac{\ln(T)}{T^{ \gamma + \rho(1-\beta) }} \right) \; \as, \qquad \text{and} \qquad  \left\| \theta_{T,w}  - \theta^{*} \right\|^{2} = \mathcal{O} \left( \frac{\ln(T)}{T^{ \gamma + \rho(1-\beta) }} \right) \; \as,
\] 
which can also be written as
\[
\left\| \theta_{T} - \theta^{*} \right\|^{2} = \mathcal{O} \left( \frac{\ln(N_{T})}{N_{T}^{\frac{\gamma + \rho(1-\beta)}{1+ \rho}}} \right) \; \as, \qquad \text{and} \qquad \left\| \theta_{T,w} - \theta^{*} \right\|^{2} = \mathcal{O} \left( \frac{\ln(N_{T})}{N_{T}^{\frac{\gamma + \rho(1-\beta)}{1+ \rho}}} \right) \; 
\as
\]

\paragraph{Rate of convergence of $\theta_{T,w}$.} Let us consider the event:
\[
E_{t} = \left\lbrace \exists j , G_{t }^{(j)} \neq \left( \frac{1}{N_{T}}\left(  G_{0 }^{(j)} + \sum_{t=1}^{T}\sum_{i=1}^{n_{t}} \left( \frac{\partial}{\partial j}f_{t,i} \left( \theta_{t-1,w} \right)\right)^{2}   \right)\right)^{-1/2}  \right\rbrace,
\]
where $f_{t,i}\left( \theta_{t-1,w} \right) :=f \left(\theta_{t-1,w};\xi_{t,i} \right)$. Observe that since $G_{t }$ converges to $G $, $\mathbb{1}_{\{E_{t}\}}$ converges almost surely to $0$, such that 
 \[
\frac{1}{s_{T}}   \sum_{t=0}^{T-1} n_{t+1} \ln (t+1)^{w + 1/2 + \delta} \left\|  G_{t+1}^{-1} - G_{t}^{-1} \right\|_{\op} \mathbb{1}_{\{E_{t}\bigcup E_{t+1}\}}  (t+1)^{\frac{\gamma -  \rho (\beta +1)}{2}} = \mathcal{O} \left( \frac{1}{T^{(1+ \rho )} \ln (T)^{w }   } \right) \; \as
 \]
 In addition, on $\{E_{t}^{C} \bigcap E_{t+1}^{C}\}$, one has 
 \begin{align*}
&\left(  G_{t+1 }^{-1} - G_{t }^{-1} \right) \mathbb{1}_{\{E_{t}^{C} \bigcap E_{t+1}^{C}\}}    =  \left( G_{t+1 }^{-1}  + G_{t }^{-1} \right)^{-1} \left( G_{t+1 }^{-2} - G_{t }^{-2} \right) \mathbb{1}_{E_{t}^{C} \bigcap E_{t+1}^{C}}  \\
& =  \left( G_{t }^{-1} + G_{t+1 }^{-1} \right)^{-1} \frac{1}{N_{t+1}}  \left( \text{diag} \left(  \sum_{i=1}^{n_{t+1}}  \left( \frac{\partial}{\partial j} f_{ t+1,i} \left( \theta_{t,w} \right) \right)^{2} \right)_{j=1,\dots ,d} - n_{t+1} G_{t }^{-2} \right) \mathbb{1}_{\{E_{t}^{C} \bigcap E_{t+1}^{C}\}},
\end{align*}
where $\text{diag} \left(  \sum_{i=1}^{n_{t+1}}  \left( \frac{\partial}{\partial j} f_{ t+1,i} \left( \theta_{t,w} \right) \right)^{2} \right)_{j=1,\dots ,d}$ is the diagonal matrix whose elements are $   \sum_{i=1}^{n_{t+1}}  \left( \frac{\partial}{\partial j} f_{ t+1,i} \left( \theta_{t,w} \right) \right)^{2}  $.
Observe that since $G_{t }$ converges almost surely to a positive matrix, there are positive constants $c_{ada},C_{ada}$ such that $\mathbb{1}_{\{E_{t,1}\}}$ converges almost surely to $1$, where $E_{t,1} := \left\lbrace  c_{ada} < \lambda_{\min} \left( G_{t } \right) \leq \lambda_{\max} \left( G_{t} \right) < C_{ada} \right\rbrace$. Then,
\begin{align*}
& \lVert (G_{t+1} - G_{t}^{-1})^{-1} \rVert_{\op} \mathbb{1}_{\{E_{t}^{C} \bigcap E_{t+1}^{C}\}} \\ & \leq  \lVert G_{t+1}^{-1} + G_{t }^{-1} \rVert_{\op} \frac{1}{N_{t+1}} \left\| \text{diag}  \left( \sum_{i=1}^{n_{t+1}}  \left( \frac{\partial}{\partial j} f_{ t+1,i} \left( \theta_{t,w} \right) \right)^{2} \right) - n_{t+1} G_{t }^{-2} \right\|_{\op} \mathbb{1}_{\{E_{t,1}^{C} \bigcup E_{t+1,1}^{C}\}} \\
& + 2C_{ada} \frac{1}{N_{t+1}} \left( \sum_{j=1}^{d} \sum_{i=1}^{n_{t+1}}  \left( \frac{\partial}{\partial j} f_{t+1,i} \left( \theta_{t,w} \right) \right)^{2} + n_{t+1}c_{ada}^{-2} \right).
\end{align*}
In addition, since $\mathbb{1}_{\{E_{t,1}^{C}\}}$ converges almost surely to $0$, one can easily check that 
\begin{align*}
&  \frac{1}{s_{T}}    \sum_{t=0}^{T-1} n_{t+1} \ln (t+1)^{w + 1/2 + \delta} (t+1)^{\frac{\gamma - ( \beta +1) \rho}{2}} \\
&  \times \left\| G_{t }^{-1} + G_{t+1 }^{-1} \right\|_{\op} \frac{1}{N_{t+1}} \left\| \text{diag}  \left( \sum_{i=1}^{n_{t+1}}  \left( \frac{\partial}{\partial j} f_{t+1,i} \left( \theta_{t,w} \right) \right)^{2}\right) - n_{t+1} G_{t }^{-2}  \right\|_{\op} \mathbb{1}_{\{E_{t,1}^{C} \bigcup E_{t+1,1}^{C}\}} \\ & = \mathcal{O} \left( \frac{1}{T^{(1+ \rho )} \ln(T)^{w }   } \right) \; \as
\end{align*}
In addition,
\[
\frac{1}{s_{T}}    \sum_{t=0}^{T-1} n_{t+1} \ln (t+1)^{w + 1/2 + \delta} (t+1)^{\frac{\gamma -( \beta+1) \rho}{2}}2C_{ada}c_{ada}^{-2} \frac{n_{t+1}}{N_{t+1}}   = O \left( \frac{\ln (T)^{\delta + 1/2 - w}}{T^\frac{2-\gamma + (\beta +1) \rho}{2}} \right) \; \as,
\]
which is negligible as soon as $ \gamma - \beta \rho < 1$. In addition, remark that 
\[
  \frac{1}{N_{t+1}} \sum_{j=1}^{d}\sum_{i=1}^{n_{t+1}}  \left(    \frac{\partial}{\partial j} f_{ t+1,i} \left( \theta_{t,w} \right) \right)^{2} =   \frac{1}{N_{t+1}} \sum_{j=1}^{d} \mathbb{E}\left[ \sum_{i=1}^{n_{t+1}}   \left(   \frac{\partial}{\partial j} f_{ t+1,i} \left( \theta_{t,w} \right) \right)^{2} 
  \middle\vert \mathcal{F}_{t} \right] +   \frac{n_{t+1}}{N_{t+1}}  \tilde{\xi}_{t+1},
\]
with $\tilde{\xi}_{t+1} =   \sum_{j=1}^{d} \left( \sum_{i=1}^{n_{t+1}}   \left(   \frac{\partial}{\partial j} f_{ t+1,i} \left( \theta_{t,w} \right) \right)^{2}   - \mathbb{E}\left[ \sum_{i=1}^{n_{t+1}}   \left(   \frac{\partial}{\partial j} f_{ t+1,i} \left( \theta_{t,w} \right) \right)^{2} \middle\vert\mathcal{F}_{t} \right] \right)$. Since $\theta_{t,w}$ converges almost surely to $\theta^{*}$ and with the help of inequality \cref{eq::momenteta}, one has
\begin{align*}
\frac{1}{s_{T}} &   \sum_{t=0}^{T-1} n_{t+1} \ln (t+1)^{w + 1/2 + \delta} (t+1)^{\frac{\gamma -( \beta+1) \rho}{2}}  \frac{1}{N_{t+1}} \sum_{j=1}^{d} \mathbb{E}\left[ \sum_{i=1}^{n_{t+1}}   \left(   \frac{\partial}{\partial j}f_{t+1,i} \left( \theta_{t,w} \right) \right)^{2} |\mathcal{F}_{t} \right]  \\ & = \mathcal{O} \left( \frac{\ln (T)^{\delta + 1/2 - w}}{T^\frac{2-\gamma + (\beta +1) \rho}{2}} \right) \; \as,
\end{align*}
while with the help of a law of large numbers for martingales (e.g., see \citet{duflo2013random}), one has
\[ 
\frac{1}{s_{T}}    \sum_{t=0}^{T-1} n_{t+1} \ln (t+1)^{w + 1/2 + \delta} (t+1)^{\frac{\gamma -( \beta+1) \rho}{2}}\frac{n_{t+1}}{N_{t+1}}  \tilde{\xi}_{t+1} = o \left( \frac{\ln (T)^{\delta + 1/2 - w}}{T^\frac{2-\gamma + (\beta +1) \rho}{2}} \right) \; \as,
\]
which concludes the proof.

\vskip 0.2in
\bibliographystyle{apalike}
\bibliography{references}

\end{document}